\documentclass{article}

\usepackage{amssymb, amsmath}
\usepackage{cite}

\begin{document}

%%% Global strong solutions of the coupled Klein-Gordon-Schr\"{o}dinger equations %%%

\begin{center}
   \LARGE \textbf{Global strong solutions of the coupled Klein-Gordon-Schr\"{o}dinger equations}
\end{center}
\vspace{5pt}
\begin{center}
   \large Tohru Ozawa and Kenta Tomioka \\
Department of Applied Physics, Waseda University   \\
Tokyo 169-8555, Japan
\end{center}
\vspace{5pt}
\begin{center}
   \textbf{Abstract.}
\end{center}

We study the initial-boundary value problem for the coupled Klein-Gordon-Schr\"{o}dinger equations in a domain in $\mathbb R^N$ with $N \leq 4$.
Under natural assumptions on the initial data, we prove the existence and uniqueness of global solutions in $H^2 \oplus H^2 \oplus H^1$.
The method of the construction of global strong solutions depends on the proof that solutions of regularized systems by the Yosida approximation form a bounded sequence in $H^2 \oplus H^2 \oplus H^1$ and a convergent sequence in $H^1 \oplus H^1 \oplus L^2$.
The method of proof is independent of the Brezis-Gallouet technique and a compactness argument.   \\



\pagebreak

%%%%%%%%%% 1. Introduction %%%%%%%%%%

\section{Introduction}

\ \ \ In this paper we study the initial-boundary value problem for the Klein-Gordon-Schr\"{o}dinger equations with Yukawa coupling in a (not necessarily bounded) domain $\Omega \subset \mathbb R^N$ with smooth boundary $\partial \Omega$ and $N \leq 4$.
In dimensionless form, the Klein-Gordon-Schr\"{o}dinger equations with Yukawa coupling are formulated as a system of equations in $\mathbb R \times \Omega$ as
\begin{align*}
   &i \partial_t u + \Delta u = - v u,   \tag{1.1}    \\
   &\partial_t^2 v - \Delta v + v = |u|^2   \tag{1.2}
\end{align*}
with initial condition at $t = 0$
\begin{align*}
   (u(0),v(0),\partial_t v(0)) = (\varphi,\psi_0,\psi_1)   \tag{1.3}
\end{align*}
and homogeneous Dirichlet boundary condition on $\partial \Omega$
\begin{align*}
   u|_{\partial \Omega} = v|_{\partial \Omega} = 0,   \tag{1.4}
\end{align*}
where $u:\mathbb R \times \Omega \to \mathbb C$ describes the complex scalar nucleon field,  $v:\mathbb R \times \Omega \to \mathbb R$ denotes the real scalar meson field, and $\varphi:\Omega \to \mathbb C,\ \psi_0:\Omega \to \mathbb R$, and $\psi_1:\Omega \to \mathbb R$ are given functions as the initial data.
We call $(u,v)$ a pair of global strong solutions to the Klein-Gordon-Schr\"{o}dinger equations if and only if $(u,v)$ satisfies (1.1) in $L^2(\Omega)$ with
\begin{align*}
   &u \in C(\mathbb R ; D(\Delta)) \cap C^1(\mathbb R ; L^2(\Omega)),   \tag{1.5}   \\
   &v \in C(\mathbb R ; D(\Delta)) \cap C^2(\mathbb R ; L^2(\Omega)).   \tag{1.6}
\end{align*}
Here the Laplacian $\Delta$ is understood to be the self-adjoint realization in the Hilbert space $L^2(\Omega)$ with domain $D(\Delta) := (H^2\cap H_0^1)(\Omega)$, where $H^2(\Omega)$ is the Sobolev space of the second order and $H_0^1(\Omega)$ is the Sobolev space of the first order satisfying the vanishing condition (1.4) on $\partial \Omega$.
By definition, strong solutions automatically satisfy (1.4).

There is a large literature on the Klein-Gordon-Schr\"{o}dinger equations in the special case where $\Omega = \mathbb R^N$ \cite{B1,CDC,CO,CHT,FT1,MX1,MX2,P1,P2,P3,R,RS,SLW,SWL,SWLW,T}, while there are few papers treating (1.1)-(1.4) in $\mathbb R \times \Omega$ with a general domain $\Omega \subset \mathbb R^N$ \cite{B2,FT2,FT3,Hn,HW}.
A major reason consists in the fact methods based on the Fourier transform are no longer available in a domain $\Omega \subsetneqq \mathbb R^N$.
Here we recall the preceding basic results on the problem (1.1)-(1.4).
In \cite{FT2,FT3}, Fukuda and Tsutsumi proved the existence and uniqueness of global strong solutions to (1.1)-(1.4) for the data
\begin{align*}
   (\varphi, \psi_0, \psi_1) \in (H^3 \cap H_0^1)(\Omega) \oplus D(\Delta) \oplus H_0^1(\Omega)
\end{align*}
in a bounded domain $\Omega \subset \mathbb R^3$.
The proof depends on the Galerkin method, where a compactness argument plays an essential role.
In \cite{HW}, Hayashi and von Wahl proved the exsistence and uniqeness of global strong solutions to (1.1)-(1.4) for the data
\begin{align*}
   (\varphi, \psi_0, \psi_1) \in D(\Delta) \oplus D(\Delta) \oplus H_0^1(\Omega)   \tag{1.7}
\end{align*}
in a (not necessarily bounded) domain $\Omega \subset \mathbb R^3$.
The proof depends on bilinear estimates in Besov spaces, nonlinear interpolation method, and the Brezis-Gallouet inequality.

The purpose of this paper is to prove the existence and uniqueness of global strong solutions to (1.1)-(1.4) for the data in the same class as in (1.7) in a (not necessarily bounded) domain $\Omega \subset \mathbb R^N$ with $N \leq 4$.
The proof depends on a modified energy argument \cite{HO1,HO2,OT1,OT2,OV} and is independent of a compactness argument, Besov estimates, nonlinear interpolation, and the Brezis-Gallouet inequality.
Moreover, we present a new information on an explicit dependence on time of the $H^2 \oplus H^2 \oplus H^1$-bounded of global strong solutions of (1.1)-(1.4).

To state our main results precisely, we introduce some notation.
The total charge and energy are defined by
\begin{align*}
   &Q(t) := \|u(t)\|_2^2,   \\
   &E(t) := \|\nabla u(t)\|_2^2 
   + \frac{1}{2} (\|\nabla v(t)\|_2^2 + \|v(t)\|_2^2 + \|\partial_t v(t)\|_2^2)
   - (v(t) | |u(t)|^2),
\end{align*}
where $(\cdot | \cdot)$ is the scalar product in $L^2$.

We now state our main theorems:   \\

\noindent
\textbf{Theorem 1.}\ \ \textit{Let $N \leq 3$ and $(\varphi,\psi_0,\psi_1) \in D(\Delta) \oplus D(\Delta) \oplus H_0^1(\Omega)$. 
Then}:   \\
\noindent
(1) \textit{There exists a unique pair of solutions $(u, v)$ satisfying} (1.1)-(1.4),
\begin{align*}
   &u \in C(\mathbb R ; D(\Delta)) \cap C^1(\mathbb R ; L^2),   \\
   &v \in C(\mathbb R ; D(\Delta)) \cap C^1(\mathbb R ; H_0^1) 
   \cap C^2(\mathbb R ; L^2).
\end{align*}
(2) \textit{The total charge $Q(t)$ and energy $E(t)$ are conserved in time.}   \\
(3) \textit{There exists a constant $C$ depending on $\max(\|\varphi\|_{H^1}, \|\psi_0\|_{H^1}, \|\psi_1\|_2)$ such that}
\begin{align*}
   \|u(t)\|_{H^1}^2 + \|v(t)\|_{H^1}^2 + \|\partial_t v(t)\|_2^2 \leq C   \tag{1.8}
\end{align*}
\textit{for all $t \in \mathbb R$.}   \\
(4) \textit{There exists a constant $C'$ depending on $\max(\|\varphi\|_{H^2}, \|\psi_0\|_{H^2}, \|\psi_1\|_{H^1})$ and a constant $C$ depending on $\max(\|\varphi\|_{H^1}, \|\psi_0\|_{H^1}, \|\psi_1\|_2)$ such that}
\begin{align*}
   \hspace{20pt} \|u(t)\|_{H^2}^2 + \|v(t)\|_{H^2}^2 + \|\partial_t v(t)\|_{H^1}^2 \leq
   \begin{cases}
   C' + C|t|^{4/3}   &\textit{for}\ N=1,   \hspace{23pt} (1.9)   \\
   C' + Ct^2   &\textit{for}\ N=2,   \hspace{18pt} (1.10)   \\
   C' + Ct^4   &\textit{for}\ N=3   \hspace{22pt} (1.11)   
   \end{cases}
\end{align*}
\textit{for all $t \in \mathbb R$.}  \\

\noindent
\textbf{Theorem 2.}\ \ \textit{Let $N = 4$.
Then there exists $\delta>0$ such that for any $(\varphi,\psi_0,\psi_1) \in D(\Delta) \oplus D(\Delta) \oplus H_0^1(\Omega)$ with $\|\varphi\|_2 < \delta$ the following holds}:   \\
\noindent
(1) \textit{There exists a unique pair of solutions $(u, v)$ satisfying} (1.1)-(1.4),
\begin{align*}
   &u \in C_w(\mathbb R ; D(\Delta)) \cap C_w^1(\mathbb R ; L^2),   \\
   &v \in C_w(\mathbb R ; D(\Delta)) \cap C_w^1(\mathbb R ; H_0^1) 
   \cap C_w^2(\mathbb R ; L^2).
\end{align*}
(2) \textit{The total charge $Q(t)$ and energy $E(t)$ are conserved in time.}   \\
(3) \textit{There exists a constant $C$ depending on $\max(\|\varphi\|_{H^1}, \|\psi_0\|_{H^1}, \|\psi_1\|_2)$ and $\delta$ such that}
\begin{align*}
   \|u(t)\|_{H^1}^2 + \|v(t)\|_{H^1}^2 + \|\partial_t v(t)\|_2^2 \leq C   \tag{1.12}
\end{align*}
\textit{for all $t \in \mathbb R$.}   \\
(4) \textit{There exist constants $C, C'$ depending on $\delta$, where $C$ and $C'$ depends on $\max(\|\varphi\|_{H^2}, \|\psi_0\|_{H^2}, \|\psi_1\|_{H^1})$ and $\max(\|\varphi\|_{H^1}, \|\psi_0\|_{H^1}, \|\psi_1\|_2)$,respectively, such that}
\begin{align*}
   \|u(t)\|_{H^2}^2 + \|v(t)\|_{H^2}^2 + \|\partial_t v(t)\|_{H^1}^2 
   \leq C \exp(C' |t|)   \tag{1.13} 
\end{align*}
\textit{for all $t \in \mathbb R$.}  \\

\noindent
\textbf{Remark 1.}\ \ (1) Boundedness of the domain $\Omega$ is not assumed in the theorems since our method is independent of a compactness argument.   \\
(2) The estimate (1.11) is an improvement of the estimates (3.13)-(3.14) in \cite{HW}, where the LHS are bounded by $C \exp(C \exp(C |t|))$ as a natural consequence of the application of the Brezis-Gallouet inequality \cite{BG,Oz}.   \\
(3) The constant $\delta$ is explicitly given by $\displaystyle \delta = \min \left( \frac{1}{C_{4,4} (C_{4,8/3})^2}, \frac{1}{\sqrt{2} C_{4,4}^2} \right)$, where $C_{N,p}$ is the best constant in the Gagliardo-Nirenberg inequalities of the form
\begin{align*}
   \|u\|_p \leq C_{N,p} \|u\|_{H^1}^{\delta_{N}(p)} \|u\|_2^{1-\delta_{N}(p)}   \tag{1.14}
\end{align*}
for all $u \in H_0^1(\Omega)$, where $\delta_{N}(p) := \frac{N}{2} - \frac{N}{p} \in [0,1]$.
See \cite{C,CH}.   \\

We prove the theorems in Section 5.
Here we illustrate the outline of the proof.
The first step of the proof is the introduction of the regularized system with $n \in \mathbb Z_{>0}$ as
\begin{align*}
   &i \partial_t u + \Delta u = - J_n^2 (J_n^2 v \cdot J_n^2 u),   \tag{1.15}   \\
   &\partial_t^2 v - \Delta v + v = J_n^2 |J_n^2 u|^2,   \tag{1.16}
\end{align*}
where $J_n := (I - \frac{1}{n} \Delta)^{-1}$ is the Yosida approximation of the identity operator $I$ given by a constant multiple of the resolvent of the Laplacian.
The associated initial data is also regularized as
\begin{align*}
   (u(0), v(0), \partial_t v(0)) = (J_n^2 \varphi, J_n^2 \psi_0, J_n^2 \psi_1).   \tag{1.17}
\end{align*}
The integural equations corresponding to (1.15)-(1.17) are given by
\begin{align*}
   &u(t) = U(t) J_n^2 \varphi 
   + i \int_0^t U(t-s) J_n^2 (J_n^2 v(s) J_n^2 u(s)) ds,   \tag{1.18}   \\
   &v(t) = \dot{K}(t) J_n^2 \psi_0 + K(t) J_n^2 \psi_1 
   + \int_0^t K(t-s) J_n^2 |J_n^2 u(s)|^2 ds,   \tag{1.19}
\end{align*}
where $U(t) := \exp(it\Delta), K(t) := \omega^{-1} \sin(t\omega), \dot{K}(t) := \cos(t\omega)$, and $\omega = (I-\Delta)^{1/2}$.
For any $n \in \mathbb Z_{>0}$, (1.18)-(1.19) have a unique local solution $(u_n, v_n) \in C(I_n ; D(\Delta^2) \oplus D(\Delta^2)) \cap C^2(I_n ; D(\Delta) \oplus D(\Delta))$, where $I_n$ is an interval containing 0 in its interior, since the regularized interaction terms are locally Lipschitz on $D(\Delta^2) \oplus D(\Delta)$.
Then a formal integration by parts based on (1.15)-(1.16) are justified on the basis of a sufficient regularity of $(u_n, v_n)$.
The second step of the proof is to obtain the $H^2 \oplus H^2 \oplus H^1$ bounded of $(u_n, v_n, \partial_t v_n)$ that is uniform in $n$ (this uniform bound ensures the global existence of strong solutions $(u_n, v_n, \partial_t v_n)$).
For that purpose we introduce a modified energy at the level of $H^2 \oplus H^2 \oplus H^1$ (see \cite{FMO,Hm,HO1,HO2,OT1,OT2,OV} for related topics).
The third step of the proof is to prove the convergence of the sequence $((u_n, v_n, \partial_t v_n) ; n \geq 1)$ in $H^1 \oplus H^1 \oplus L^2$.
The proof that the sequence $((u_n, v_n, \partial_t v_n) ; n \geq 1)$ forms a Cauchy sequence depends essentially on the uniform boundedness in $H^2 \oplus H^2 \oplus H^1$.
The last step of the proof is the existence and uniqueness of global strong solutions to (1.1)-(1.4) satisfying the statements of the theorems.

We close the introduction by explaining the organization of this paper.
In Section 2, we prepare basic lemmas for the proofs of the main theorems.
In Section 3, we prove the uniform boundedness of the sequence of solutions $((u_n, v_n, \partial_t v_n) ; n \geq 1)$ to (1.18)-(1.19) in $H^2 \oplus H^2 \oplus H^1$.
In Section 4, we prove the convergence of the sequence $((u_n, v_n, \partial_t v_n) ; n \geq 1)$ in $H^1 \oplus H^1 \oplus L^2$.
In Section 5, we prove Theorems 1 and 2.
In Section 6, we give remarks on the existence and uniqueness of global finite energy solutions.

%%%%%%%%%% 2. Preliminaries %%%%%%%%%%   

\section{Preliminaries}

\ \ \ In this section we collect basic estimates and propositions for the proofs of the main theorems.
Different positive constants independent of the time variable $t$ might be denoted by the same letter $C$.   \\

\noindent
\textbf{Lemma 1}(Elliptic estimate \cite{C,CH})\textbf{.}\ \ \textit{There exists a constant $C$ such that the estimate}
\begin{align*}
   &\|u\|_{H^2} \leq C \|(I - \Delta) u\|_2   \tag{2.1}
\end{align*}
holds for any $u \in D(\Delta)$.   \\

\noindent
\textbf{Lemma 2}(Gagliardo-Nirenberg inequalities \cite{C,CH})\textbf{.}\ \ \textit{For any $N$ and $p$ with $\delta_{N}(p) := \frac{N}{2} - \frac{N}{p} \in [0,1]$, there exists a constant $C_{N,p}$ such that the estimate}
\begin{align*}
   \|u\|_p \leq C_{N,p} \|u\|_{H^1}^{\delta_{N}(p)} \|u\|_2^{1 - \delta_{N}(p)}   \tag{2.2}
\end{align*}
holds for any $u \in H_0^1(\Omega)$.   \\

\noindent
\textbf{Lemma 3}(\cite{C})\textbf{.}\ \ \textit{The operators $J_n = (I- \frac{1}{n} \Delta)^{-1}$ satisfy the following properties}: \\
\noindent
(1) \textit{$J_n$ are bounded self-adjoint operators in $L^2(\Omega)$ with image}
\begin{align*}
   J_n (L^2(\Omega)) = D(\Delta).   
\end{align*}
(2) \textit{For any $u\in L^2(\Omega), J_n u \to u$ in $L^2(\Omega)$ as $n \to \infty$.}   \\
(3) \textit{For any $u \in L^2(\Omega)$, the following estimates hold.}
\begin{align*}
   &\|J_n u\|_2 \leq \|u\|_2,   \tag{2.3}   \\[-2pt]
   &\|\nabla J_n u\|_2 \leq n^{\frac{1}{2}} \|u\|_2,  \tag{2.4}   \\
   &\|\Delta J_n u\|_2 \leq n \|u\|_2.   \tag{2.5}   
\end{align*}
(4) \textit{For any $m, n$ with $m>n$ and any $u \in H_0^1(\Omega)$, the following estimate holds.}
\begin{align*}
   \|(J_m - J_n) u\|_2 \leq \frac{1}{n^{\frac{1}{2}}} \|\nabla u\|_2.   \tag{2.6}
\end{align*}

%% 3. Uniform boundedness of the sequence of solutions of the regularized system %%%

\section{Uniform boundedness of the sequence of \\ 
solutions of the regularized system}

\ \ \ Let $(\varphi, \psi_0, \psi_1) \in D(\Delta) \oplus D(\Delta) \oplus H_0^1(\Omega)$.
We assume \\
$\displaystyle \|\varphi\|_2 < \min \left( \frac{1}{C_{4,4} (C_{4,3/8})^2}, \frac{1}{\sqrt{2} C_{4,4}^2} \right)$ for $N=4$.
For any $n \in \mathbb Z_{>0}$, (1.18)-(1.19) have a unique local solution $(u_n, v_n) \in C(I_n ; D(\Delta^2) \oplus D(\Delta^2)) \cap C^2(I_n ; D(\Delta) \oplus D(\Delta))$, where $I_n$ is an interval containing 0 in its interior.
Then $(u_n, v_n)$ satisfies
\begin{align*}
   &i \partial_t u_n + \Delta u_n = - J_n^2 (J_n^2 v_n \cdot J_n^2 u_n),   \tag{3.1}    \\
   &\partial_t^2 v_n - \Delta v_n + v_n = J_n^2 |J_n^2 u_n|^2   \tag{3.2}   \\
   &(u_n(0), v_n(0), \partial_t v_n(0)) = (J_n^2 \varphi, J_n^2 \psi_0, J_n^2 \psi_1).   \tag{3.3}
\end{align*}
Omitting explicit dependence of the time variable as usual, by (3.1) and (3.2) we have
\begin{align*}
   \frac{d}{dt} \|u_n\|_2^2 
   &= 2 \textrm{Re} (\partial_t u_n | u_n) = 2 \textrm{Im} (i \partial_t u_n | u_n)   \\[-5pt]
   &= 2 \textrm{Im} (-\Delta u_n - J_n^2 (J_n^2 v_n \cdot J_n^2 u_n) | u_n)   \\
   &= 2 \textrm{Im} \|\nabla u_n\|_2^2 
   - 2 \textrm{Im} (J_n^2 v_n \cdot J_n^2 u_n | J_n^2 u_n) = 0,   \tag{3.4}
\end{align*}
\begin{align*}
   &\frac{d}{dt} (\|\nabla u_n\|_2^2 
   +\frac{1}{2} (\|\partial_t v_n\|_2^2 + \|\nabla v_n \|_2^2 + \|v_n\|_2^2))    \\[-2pt]
   &= 2 \textrm{Re} (\partial_t \nabla u_n | \nabla u_n) + (\partial_t^2 v_n | \partial_t v_n )
   + (\partial_t \nabla v_n | \nabla v_n) + (\partial_t v_n | v_n)   \\
   &= - 2 \textrm{Re} (\partial_t u_n | \Delta u_n) 
   + (\partial_t^2 v_n - \Delta v_n + v_n | \partial_t v_n)  \\
   &= - 2 \textrm{Re} (\partial_t u_n | - i \partial_t u_n - J_n^2 (J_n^2 v_n \cdot J_n^2 u_n))
   + (J_n^2 |J_n^2 u_n|^2 | \partial_t v_n)   \\
   &= (\partial_t |J_n^2 u_n|^2 | J_n^2 v_n) + (|J_n^2 u_n|^2 | \partial_t J_n^2 v_n)   \\
   &= \frac{d}{dt} (|J_n^2 u_n|^2 | J_n^2 v_n),   \tag{3.5}
\end{align*}
from which we obtain the conservation laws:
\begin{align*}
   &\|u_n (t)\|_2^2 = \|u_n (0)\|_2^2 = \|J_n^2 \varphi \|_2^2,   \tag{3.6}   \\
   &E_n (t) = E_n (0),   \tag{3.7}
\end{align*}
where
\begin{align*}
   E_n (t) := \|\nabla u_n\|_2^2 
   + \frac{1}{2} (\|\partial_t v_n\|_2^2 + \|\nabla v_n\|_2^2 + \|v_n\|_2^2) 
   - (|J_n^2 u_n|^2 | J_n^2 v_n).   \tag{3.8}
\end{align*}
By (3.6), we have
\begin{align*}
   \|u_n(t)\|_2^2 \leq \|\varphi\|_2^2.   \tag{3.9}
\end{align*}
By the H\"{o}lder and Gagliardo-Nirenberg inequalities (Lemma 2), Lemma 3, and (3.9), the last term on the RHS of (3.8) is estimated by
\begin{align*}
   |(|J_n^2 u_n |^2 | J_n^2 v_n)| 
   &\leq \|J_n^2 u_n \|_2^2 \|J_n^2 v_n \|_{\infty}   \\
   &\leq C \|u_n \|_2^2 \|v_n \|_{H^1}   \\
   &\leq C \|\varphi \|_2^2 (\|\nabla v_n \|_2^2 + \|v_n\|_2^2)^{1/2}   \\
   &\leq \frac{1}{4} (\|\nabla v_n\|_2^2 + \|v_n\|_2^2) + C \|\varphi \|_2^4
   \hspace{25pt}  \textrm{for}\ N=1,   \tag{3.10}   \\
   |(|J_n^2 u_n |^2 | J_n^2 v_n)| 
   &\leq \|J_n^2 u_n \|_{8/3}^2 \|J_n^2 v_n \|_4   \\
   &\leq C \|u_n \|_{H^1}^{1/2} \|u_n\|_2^{3/2} \|v_n \|_{H^1}   \\
   &\leq C (\|\nabla u_n\|_2^2 + \|\varphi\|_2^2)^{1/4} \|\varphi \|_2^{3/2}
   (\|\nabla v_n \|_2^2 + \|v_n\|_2^2)^{1/2}   \\
   &\leq \frac{1}{2} \|\nabla u_n\|_2^2 
   + \frac{1}{4} (\|\nabla v_n\|_2^2 + \|v_n\|_2^2) 
   + C (\|\varphi\|_2^2 + \|\varphi \|_2^6)   \\
   &\hspace{165pt} \textrm{for}\ N=2,   \tag{3.11}
\end{align*}
\begin{align*}
   \hspace{13pt} |(|J_n^2 u_n |^2 | J_n^2 v_n)| 
   &\leq \|J_n^2 u_n \|_{12/5}^2 \|J_n^2 v_n \|_6   \\
   &\leq C \|u_n \|_{H^1}^{1/2} \|u_n\|_2^{3/2} \|v_n \|_{H^1}   \\
   &\leq C (\|\nabla u_n\|_2^2 + \|\varphi\|_2^2)^{1/4} \|\varphi \|_2^{3/2}
   (\|\nabla v_n \|_2^2 + \|v_n\|_2^2)^{1/2}   \\
   &\leq \frac{1}{2} \|\nabla u_n\|_2^2 
   + \frac{1}{4} (\|\nabla v_n\|_2^2 + \|v_n\|_2^2) 
   + C (\|\varphi\|_2^2 + \|\varphi \|_2^6)   \\
   &\hspace{165pt} \textrm{for}\ N=3,   \tag{3.12}   \\
   |(|J_n^2 u_n |^2 | J_n^2 v_n)| 
   &\leq \|J_n^2 u_n \|_{8/3}^2 \|J_n^2 v_n \|_4   \\
   &\leq (C_{4,8/3}^2 \|u_n \|_{H^1} \|u_n\|_2) (C_{4,4} \|v_n \|_{H^1})   \\
   &\leq \frac{1}{4} (\|\nabla v_n\|_2^2 + \|v_n\|_2^2) 
   + C_{4,4}^2 C_{4,8/3}^4 (\|\nabla u_n\|_2^2 + \|\varphi \|_2^2) \|\varphi\|_2^2   \\
   &\hspace{165pt} \textrm{for}\ N=4,   \tag{3.13}
\end{align*}
By (3.10)-(3.13), $E_n(t)$ is bounded from below by
\begin{align*}
   \hspace{6pt} E_n(t) \geq
   \begin{cases}
   \frac{1}{2} \|\nabla u_n\|_2^2 
   + \frac{1}{4} (\|\partial_t v_n\|_2^2 + \|\nabla v_n\|_2^2 + \|v_n\|_2^2)
   - C \|\varphi\|_2^4   \\
   \hspace{200pt} \textrm{for}\ N=1, \hspace{20pt} (3.14)   \\
   \frac{1}{2} \|\nabla u_n\|_2^2 
   + \frac{1}{4} (\|\partial_t v_n\|_2^2 + \|\nabla v_n\|_2^2 + \|v_n\|_2^2)
   - C (\|\varphi\|_2^2 + \|\varphi\|_2^6)   \\
   \hspace{200pt} \textrm{for}\ N=2,3, \hspace{11pt} (3.15)   \\
   (1 - C_{4,4}^2 C_{4,8/3}^4 \|\varphi\|_2^2) \|\nabla u_n\|_2^2   \\
   \hspace{11pt} + \frac{1}{4} (\|\partial_t v_n\|_2^2 + \|\nabla v_n\|_2^2 + \|v_n\|_2^2)
   - C \|\varphi\|_4^4
   \hspace{10pt} \textrm{for}\ N=4, \hspace{20pt} (3.16)   \\
   \end{cases}
\end{align*}
while $E_n(0)$ is estimated by
\begin{align*}
   E_n(0) 
   &= \|\nabla J_n^2 \varphi\|_2^2 
   + \frac{1}{2} (\|J_n^2 \psi_1\|_2^2 + \|\nabla J_n^2 \psi_0\|_2^2 + \|J_n^2 \psi_0\|_2^2)
   - (|J_n^4 \varphi|^2 | J_n^4 \psi_0)   \\[-3pt]
   &\leq \|\nabla \varphi\|_2^2
   + \frac{1}{2} (\|\nabla \psi_0\|_2^2 + \|\psi_0\|_2^2 + \|\psi_1\|_2^2)
   + C \|\varphi\|_{H^1}^2 \|\psi_0\|_{H^1}.   \tag{3.17}
\end{align*}
By (3.7),(3.14)-(3.17), we obtain the estimate
\begin{align*}
   &\frac{1}{2} \|\nabla u_n\|_2^2 
   + \frac{1}{4} (\|\partial_t v_n\|_2^2 + \|\nabla v_n\|_2^2 + \|v_n\|_2^2)   \\
   &\leq E_n(t) + C (\|\varphi\|_2^2 + \|\varphi\|_2^6)   \\
   &= E_n(0) + C (\|\varphi\|_2^2 + \|\varphi\|_2^6)   \\
   &\leq \|\nabla \varphi\|_2^2
   + \frac{1}{2} (\|\nabla \psi_0\|_2^2 + \|\psi_0\|_2^2 + \|\psi_1\|_2^2)
   + C \|\varphi\|_{H^1}^2 \|\psi_0\|_{H^1}
   + C (\|\varphi\|_2^2 + \|\varphi\|_2^6)   \\
   &\hspace{240pt} \textrm{for}\ N=1,2,3,   \tag{3.18}
\end{align*}
\vspace{-20pt}
\begin{align*}
   &(1 - C_{4,4}^2 C_{4,8/3}^4 \|\varphi\|_2^2) \|\nabla u_n\|_2^2 
   + \frac{1}{4} (\|\partial_t v_n\|_2^2 + \|\nabla v_n\|_2^2 + \|v_n\|_2^2)   \\
   &\leq \|\nabla \varphi\|_2^2
   + \frac{1}{2} (\|\nabla \psi_0\|_2^2 + \|\psi_0\|_2^2 + \|\psi_1\|_2^2)
   + C \|\varphi\|_{H^1}^2 \|\psi_0\|_{H^1}
   + C \|\varphi\|_4^4   \\
   &\hspace{225pt} \textrm{for}\ N=4.   \tag{3.19}
\end{align*}
The estimates (3.9),(3.18), and (3.19) are uniform in both $t \in \mathbb R$ and $n \geq 1$.
This implies the global existence of finite energy solutions $(u_n, v_n, \partial_t v_n)$ as well.
We introduce the resulting $H^1 \oplus H^1 \oplus L^2$ bound as $M_1$ depending on $\max(\|\varphi\|_{H^1}, \|\psi_0\|_{H^1}, \|\psi_1\|_2)$:
\begin{align*}
   \sup_{n \geq 1} \sup_{t \in \mathbb R} 
   (\|u_n(t)\|_{H^1}^2 + \|v_n(t)\|_{H^1}^2 + \|\partial_t v_n(t)\|_2^2) \leq M_1 < \infty.
   \tag{3.20}
\end{align*}

We new consider the uniform estimate in $H^2 \oplus H^2 \oplus H^1$.
For that purpose we compute
\begin{align*}
   &\frac{d}{dt} ( \|\partial_t u_n\|_2^2 + \frac{1}{2} ( \|\nabla \partial_t v_n\|_2^2 
   + \|\Delta v_n\|_2^2 + \|\nabla v_n\|_2^2) )   \\
   &= 2 \textrm{Im} ( i \partial_t^2 u_n|\partial_t u_n) 
   + (\nabla \partial_t^2 v_n|\nabla \partial_t v_n) 
   + (\Delta \partial_t v_n|\Delta v_n)
   + (\nabla \partial_t v_n|\nabla v_n)   \\
   &= 2 \textrm{Im} 
   ( \partial_t ( - \Delta u_n - J_n^2 ( J_n^2 v_n \cdot J_n^2 u_n ) ) |\partial_t u_n)
   + ( \nabla ( \partial_t^2 v_n - \Delta v_n + v_n ) | \nabla \partial_t v_n )   \\
   &= - 2 \textrm{Im} ( J_n^2 ( J_n^2 \partial_t v_n \cdot J_n^2 u_n)) | \partial_t u_n)
   + ( \nabla J_n^2 |J_n^2 u_n|^2 | \nabla \partial_t v_n)     \\
   &= - 2 \textrm{Re} ( J_n^2 ( J_n^2 \partial_t v_n \cdot J_n^2 u_n)) | i \partial_t u_n)
   + ( \nabla |J_n^2 u_n|^2 | \nabla J_n^2 \partial_t v_n)   \\
   &= -2 \textrm{Re} ( J_n^2 ( J_n^2 \partial_t v_n \cdot J_n^2 u_n ) ) | - \Delta u_n - J_n^2 ( J_n^2 v_n \cdot J_n^2 u_n ) )
   + ( \nabla |J_n^2 u_n|^2 | \nabla J_n^2 \partial_t v_n )   \\
   &= ( J_n^2 \partial_t v_n | 2 \textrm{Re} \overline{J_n^2 u_n} \Delta J_n^2 u_n )
   + 2 \textrm{Re} ( J_n^2 ( J_n^2 \partial_t v_n \cdot J_n^2 u_n ) | J_n^2 ( J_n^2 v_n \cdot J_n^2 u_n ) )   \\
   &\hspace{11pt} + ( \nabla |J_n^2 u_n|^2 | \nabla J_n^2 \partial_t v_n)   \\
   &= ( J_n^2 \partial_t v_n | \Delta |J_n^2 u_n|^2 - 2 |\nabla J_n^2 u_n|^2 )
   + 2 \textrm{Re} ( J_n^2 ( J_n^2 \partial_t v_n \cdot J_n^2 u_n ) | J_n^2 ( J_n^2 v_n \cdot J_n^2 u_n ) )   \\
   &\hspace{11pt} - ( \Delta |J_n^2 u_n|^2 | J_n^2 \partial_t v_n)   \\
   &= - 2 ( J_n^2 \partial_t v_n | |\nabla J_n^2 u_n|^2 )
   + 2 \textrm{Re} ( J_n^2 ( J_n^2 \partial_t v_n \cdot J_n^2 u_n ) | J_n^2 ( J_n^2 v_n \cdot J_n^2 u_n ) ).   \tag{3.21}
\end{align*}
Based on the equalities
\begin{align*}
   \|\partial_t u_n\|_2^2 
   &= \| \Delta u_n + J_n^2 (J_n^2 v_n \cdot J_n^2 u_n)\|_2^2   \\
   &= \|\Delta u_n\|_2^2 
   + 2 \textrm{Re} (\Delta J_n^2 u_n | J_n^2 v_n \cdot J_n^2 u_n)
   + \|J_n^2 (J_n^2 v_n \cdot J_n^2 u_n)\|_2^2   \\
   &= \|\Delta u_n\|_2^2 
   - 2 (\nabla J_n^2 v_n | \nabla |J_n^2 u_n|^2)
   - 2 (J_n^2 v_n | |\nabla J_n^2 u_n|^2)   \\
   &\hspace{11pt} + \|J_n^2 (J_n^2 v_n \cdot J_n^2 u_n)\|_2^2,   \tag{3.22}
\end{align*}
we introduce a natural second order energy $F_n(t)$ for (3.1)-(3.2) by
\begin{align*}
   F_n(t) &:= \|\Delta u_n\|_2^2 
   + \frac{1}{2} ( \|\nabla \partial_t v_n\|_2^2 + \|\Delta v_n\|_2^2 + \|\nabla v_n\|_2^2)
   + \|J_n^2 (J_n^2 v_n \cdot J_n^2 u_n)\|_2^2   \\
   &\hspace{11pt} - 2 (\nabla J_n^2 v_n | \nabla |J_n^2 u_n|^2)
   - 2 (J_n^2 v_n | |\nabla J_n^2 u_n|^2).   \tag{3.23}
\end{align*}
Then, (3.21) is rewritten as
\begin{align*}
   F_n'(t) = - 2 ( J_n^2 \partial_t v_n | |\nabla J_n^2 u_n|^2 )
   + 2 \textrm{Re} ( J_n^2 ( J_n^2 \partial_t v_n \cdot J_n^2 u_n ) | J_n^2 ( J_n^2 v_n \cdot J_n^2 u_n ) ).   \tag{3.24}
\end{align*}
Indefinite terms on the RHS of (3.23) are estimated by
\begin{align*}
   |(\nabla J_n^2 v_n | \nabla |J_n^2 u_n|^2)|
   &\leq 2 \|\nabla v_n\|_2 \|J_n^2 u_n\|_{\infty} \|\nabla u_n\|_2   \\
   &\leq C \|\nabla v_n\|_2 \|u_n\|_{H^1}^2 \leq C M_1^{3/2},   \tag{3.25}   \\
   |(J_n^2 v_n | |\nabla J_n^2 u_n|^2)|
   &\leq \|J_n^2 v_n\|_{\infty} \|\nabla u_n\|_2^2   \\
   &\leq C \|v_n\|_{H^1}^{1/2} \|v_n\|_2^{1/2} \|\nabla u_n\|_2^2 \leq C M_1^{3/2}
   \ \ \textrm{for}\ N=1,   \tag{3.26}   \\
   |(\nabla J_n^2 v_n | \nabla |J_n^2 u_n|^2)|
   &\leq 2 \|\nabla v_n\|_2 \|\nabla J_n^2 u_n\|_4 \|J_n^2 u_n\|_4   \\
   &\leq C \|\nabla v_n\|_2 (\|u_n\|_{H^2}^{1/2} \|u_n\|_{H^1}^{1/2}) 
   (\|u_n\|_{H^1}^{1/2} \|u_n\|_2^{1/2})   \\
   &\leq C \|\nabla v_n\|_2 (\|\Delta u_n\|_2^2 + \|\varphi\|_2^2)^{1/4} 
   \|\varphi\|_2^{1/2} \|u_n\|_{H^1}   \\
   &\leq C M_1 \|\varphi\|_2^{1/2} (\|\Delta u_n\|_2^2 + \|\varphi\|_2^2)^{1/4}   \\
   &\leq \frac{1}{8} (\|\Delta u_n\|_2^2 + \|\varphi\|_2^2) 
   + C M_1^{4/3} \|\varphi\|_2^{2/3},   \tag{3.27}   \\
   |(J_n^2 v_n | |\nabla J_n^2 u_n|^2)|
   &\leq \|J_n^2 v_n\|_4 \|\nabla u_n\|_2 \|\nabla J_n^2 u_n\|_4   \\
   &\leq C \|v_n\|_{H^1} \|u_n\|_{H^1}^{3/2} \|u_n\|_{H^2}^{1/2}   \\
   &\leq C M_1^{5/4} (\|\Delta u_n\|_2^2 + \|\varphi\|_2^2)^{1/4}   \\
   &\leq \frac{1}{8} (\|\Delta u_n\|_2^2 + \|\varphi\|_2^2) + C M_1^{5/3}
   \hspace{32pt} \textrm{for}\ N=2,   \tag{3.28}   
\end{align*}
\begin{align*}
   |(\nabla J_n^2 v_n | \nabla |J_n^2 u_n|^2)|
   &\leq 2 \|\nabla v_n\|_2 \|\nabla J_n^2 u_n\|_3 \|J_n^2 u_n\|_6   \\
   &\leq C \|\nabla v_n\|_2 (\|u_n\|_{H^2}^{1/2} \|\nabla u_n\|_2^{1/2}) \|u_n\|_{H^1}   \\
   &\leq C M_1^{5/4} (\|\Delta u_n\|_2^2 + \|\varphi\|_2^2)^{1/4}   \\
   &\leq \frac{1}{8} (\|\Delta u_n\|_2^2 + \|\varphi\|_2^2) + C M_1^{5/3},   \tag{3.29}   \\
   |(J_n^2 v_n | |\nabla J_n^2 u_n|^2)|
   &\leq \|J_n^2 v_n\|_6 \|\nabla J_n^2 u_n\|_{12/5}^2   \\
   &\leq C \|v_n\|_{H^1} \|u_n\|_{H^2}^{1/2} \|u_n\|_{H^1}^{3/2}   \\
   &\leq C M_1^{5/4} (\|\Delta u_n\|_2^2 + \|\varphi\|_2^2)^{1/4}   \\
   &\leq \frac{1}{8} (\|\Delta u_n\|_2^2 + \|\varphi\|_2^2) + C M_1^{5/3}
   \hspace{32pt} \textrm{for}\ N=3,   \tag{3.30}   \\
   |(\nabla J_n^2 v_n | \nabla |J_n^2 u_n|^2)|
   &\leq 2 \|\nabla v_n\|_2 \|\nabla J_n^2 u_n\|_4 \|J_n^2 u_n\|_4   \\
   &\leq C \|\nabla v_n\|_2 \|u_n\|_{H^2} \|u_n\|_{H^1}   \\
   &\leq C M_1 (\|\Delta u_n\|_2^2 + \|\varphi\|_2^2)^{1/2}   \\
   &\leq \frac{1}{8} (\|\Delta u_n\|_2^2 + \|\varphi\|_2^2) + C M_1^2,   \tag{3.31}   \\
   |(J_n^2 v_n | |\nabla J_n^2 u_n|^2)|
   &\leq \|J_n^2 v_n\|_4 \|\nabla J_n^2 u_n\|_{8/3}^2   \\
   &\leq C \|v_n\|_{H^1} \|u_n\|_{H^2} \|u_n\|_{H^1}   \\
   &\leq C M_1 (\|\Delta u_n\|_2^2 + \|\varphi\|_2^2)^{1/2}   \\
   &\leq \frac{1}{8} (\|\Delta u_n\|_2^2 + \|\varphi\|_2^2) + C M_1^2
   \hspace{40pt} \textrm{for}\ N=4.   \tag{3.32}   
\end{align*}
By (3.23), (3.25)-(3.32), we obtain
\begin{align*}
   &\|\Delta u_n\|_2^2 
   + \frac{1}{2} ( \|\nabla \partial_t v_n\|_2^2 + \|\Delta v_n\|_2^2 + \|\nabla v_n\|_2^2)
   + \|J_n^2 (J_n^2 v_n \cdot J_n^2 u_n)\|_2^2   \\[-2pt]
   &\leq 2 F_n(t) + C (M_1^{3/2} + M_1^2)   \\[-2pt]
   &\leq 3 \|\Delta u_n\|_2^2 
   + ( \|\nabla \partial_t v_n\|_2^2 + \|\Delta v_n\|_2^2 + \|\nabla v_n\|_2^2)
   + C (M_1^{3/2} + M_1^3),   \tag{3.33}
\end{align*}
where we have used
\begin{align*}
   \hspace{20pt} \|J_n^2 (J_n^2 v_n \cdot J_n^2 u_n)\|_2^2
   &\leq \|J_n^2 v_n \cdot J_n^2 u_n\|_2^2   \\[-1pt]
   &\leq C \|v_n\|_2^2 \|J_n^2 u_n\|_{\infty}^2   \\[-1pt]
   &\leq C M_1 \|J_n^2 u_n\|_{\infty}^2   \\
   &\leq C M_1 \|u_n\|_{H^2} \|u_n\|_{H^1}   \\
   &\leq \frac{1}{4} (\|\Delta u_n\|_2^2 + \|\varphi\|_2^2) + C M_1^3
   \hspace{15pt} \textrm{for}\ N \leq 3,   \tag{3.34}   \\
   \|J_n^2 (J_n^2 v_n \cdot J_n^2 u_n)\|_2^2
   &\leq \|J_n^2 v_n\|_4^2 \|J_n^2 u_n\|_4^2   \\[-1pt]
   &\leq C \|v_n\|_{H^1}^2 \|u_n\|_{H^1}^2   \\
   &\leq C M_1^2
   \hspace{110pt} \textrm{for}\ N=4.   \tag{3.35}
\end{align*}
We estimate two terms on the RHS of (3.24) as
\begin{align*}
   |(J_n^2 \partial_t v_n | |\nabla J_n^2 u_n|^2)|
   &\leq \|J_n^2 \partial_t v_n\|_2 \|\nabla J_n^2 u_n\|_4^2   \\
   &\leq C \|\partial_t v_n\|_2 \|u_n\|_{H^2}^{1/2} \|\nabla u_n\|_2^{3/2}   \\
   &\leq C M_1^{5/4} (\|\Delta u_n\|_2^2 + \|\varphi\|_2^2)^{1/4},   \tag{3.36}   \\
   |(J_n^2 (J_n^2 \partial_t v_n \cdot J_n^2 u_n) | J_n^2 (J_n^2 v_n \cdot J_n^2 u_n))| 
   &\leq \|\partial_t v_n\|_2 \|J_n^2 u_n\|_4^2 \|J_n^2 v_n\|_{\infty}   \\
   &\leq C \|\partial_t v_n\|_2 \|u_n\|_{H^1}^{1/2} \|u_n\|_2^{3/2} \|v_n\|_{H^1}   \\
   &\leq C M_1^2 \hspace{50pt} \textrm{for}\ N=1,   \tag{3.37}
\end{align*}
\vspace{-14pt}
\begin{align*}
   |(J_n^2 \partial_t v_n | |\nabla J_n^2 u_n|^2)|
   &\leq \|\partial_t v_n\|_2 \|\nabla J_n^2 u_n\|_4^2   \\
   &\leq C \|\partial_t v_n\|_2 \|u_n\|_{H^2} \|\nabla u_n\|_2   \\
   &\leq C M_1 (\|\Delta u_n\|_2^2 + \|\varphi\|_2^2)^{1/2},   \tag{3.38}   \\
   |(J_n^2 (J_n^2 \partial_t v_n \cdot J_n^2 u_n) | J_n^2 (J_n^2 v_n \cdot J_n^2 u_n))| 
   &\leq \|\partial_t v_n\|_2 \|J_n^2 v_n\|_4 \|J_n^2 u_n\|_8^2   \\
   &\leq C \|\partial_t v_n\|_2 \|v_n\|_{H^1} \|u_n\|_{H^1}^2   \\
   &\leq C M_1^2 \hspace{50pt} \textrm{for}\ N=2, \hspace{2pt}   \tag{3.39}   
\end{align*}
\vspace{-14pt}
\begin{align*}
   |(J_n^2 \partial_t v_n | |\nabla J_n^2 u_n|^2)|
   &\leq \|J_n^2 \partial_t v_n\|_2 \|\nabla J_n^2 u_n\|_4^2   \\
   &\leq C \|\partial_t v_n\|_2 \|u_n\|_{H^2}^{3/2} \|\nabla u_n\|_2^{1/2}   \\
   &\leq C M_1^{3/4} (\|\Delta u_n\|_2^2 + \|\varphi\|_2^2)^{3/4},   \tag{3.40}   \\
   |(J_n^2 (J_n^2 \partial_t v_n \cdot J_n^2 u_n) | J_n^2 (J_n^2 v_n \cdot J_n^2 u_n))| 
   &\leq \|\partial_t v_n\|_2 \|J_n^2 v_n\|_6 \|J_n^2 u_n\|_6^2   \\
   &\leq C \|\partial_t v_n\|_2 \|v_n\|_{H^1} \|u_n\|_{H^1}^2   \\
   &\leq C M_1^2 \hspace{50pt} \textrm{for}\ N=3,   \tag{3.41}
\end{align*}
\begin{align*}
   |(J_n^2 \partial_t v_n | |\nabla J_n^2 u_n|^2)|
   &\leq \|\partial_t v_n\|_2 \|\nabla J_n^2 u_n\|_4^2   \\
   &\leq C \|\partial_t v_n\|_2 \|u_n\|_{H^2}^2   \\
   &\leq C M_1^{1/2} (\|\Delta u_n\|_2^2 + \|\varphi\|_2^2),   \tag{3.42}   \\
   |(J_n^2 (J_n^2 \partial_t v_n \cdot J_n^2 u_n) | J_n^2 (J_n^2 v_n \cdot J_n^2 u_n))| 
   &\leq \|J_n^2 \partial_t v_n\|_4 \|J_n^2 v_n\|_4 \|J_n^2 u_n\|_4^2   \\
   &\leq C \|\partial_t v_n\|_{H^1} \|v_n\|_{H^1} \|u_n\|_{H^1}^2   \\
   &\leq C (\|\nabla \partial_t v_n\|_2^2 + M_1)^{1/2} M_1^{3/2}   \\
   &\hspace{87pt} \textrm{for}\ N=4. \hspace{3pt}   \tag{3.43}
\end{align*}
By (3.24), (3.33), and (3.36)-(3.43), we obtain
\begin{align*}
   \hspace{45pt} |F_n'(t)| \leq 
   \begin{cases}
   C M_1^{5/4} (F_n(t) + M_1 + M_1^3)^{1/4} 
   \hspace{36pt} \textrm{for}\ N=1, \hspace{11pt}  (3.44)    \\
   C M_1 (F_n(t) + M_1 + M_1^3)^{1/2} 
   \hspace{45pt} \textrm{for}\ N=2, \hspace{11pt} (3.45)    \\
   C M_1^{3/4} (F_n(t) + M_1 + M_1^3)^{3/4} 
   \hspace{36pt} \textrm{for}\ N=3, \hspace{11pt} (3.46)    \\
   C M_1^{1/2} (F_n(t) + M_1 + M_1^3) 
   \hspace{49pt} \textrm{for}\ N=4. \hspace{12pt} (3.47)
   \end{cases}
\end{align*}
Replacing $F_n'(t)$ on the LHS of (3.44)-(3.47) by $(F_n + M_1 + M_1^3)'(t)$ and integrating the corresponding differential inequalities, we obtain
\begin{align*}
   \hspace{6pt} F_n(t) + M_1 + M_1^3 \leq 
   \begin{cases}
   C (F_n(0) + M_1 + M_1^3 + M_1^{5/3} |t|^{4/3}) 
   \hspace{12pt} \textrm{for}\ N=1, \hspace{4pt} (3.48)    \\
   C (F_n(0) + M_1 + M_1^3 + M_1^2 |t|^2) 
   \hspace{28pt} \textrm{for}\ N=2, \hspace{4pt} (3.49)   \\
   C (F_n(0) + M_1 + M_1^3 + M_1^3 |t|^4) 
   \hspace{28pt} \textrm{for}\ N=3, \hspace{4pt} (3.50)   \\
   C (F_n(0) + M_1 + M_1^3) \exp (C M_1^{1/2} |t|) 
   \hspace{4pt} \textrm{for}\ N=4. \hspace{5pt} (3.51)
   \end{cases}
\end{align*}
By (3.33) and (3.48)-(3.51), we obtain the estimate
\begin{align*}
   \hspace{40pt} &\|\Delta u_n\|_2^2 
   + \frac{1}{2} (\|\nabla \partial_t v_n\|_2^2 + \|\Delta v_n\|_2^2 + \|\nabla v_n\|_2^2) \\
   &\leq 
   \begin{cases}
   C (\|\Delta \varphi\|_2^2 + \|\nabla \psi_1\|_2^2 + \|\Delta \psi_0\|_2^2 
   + \|\nabla \psi_0\|_2^2   \\
   \hspace{15pt} + M_1 + M_1^3 + M_1^{5/3} |t|^{4/3})   
   \hspace{82pt} \textrm{for}\ N=1, \hspace{3pt} (3.52)   \\
   C (\|\Delta \varphi\|_2^2 + \|\nabla \psi_1\|_2^2 + \|\Delta \psi_0\|_2^2 
   + \|\nabla \psi_0\|_2^2   \\
   \hspace{15pt} + M_1 + M_1^3 + M_1^2 |t|^2)   
   \hspace{98pt} \textrm{for}\ N=2, \hspace{3pt} (3.53)   \\
   C (\|\Delta \varphi\|_2^2 + \|\nabla \psi_1\|_2^2 + \|\Delta \psi_0\|_2^2 
   + \|\nabla \psi_0\|_2^2   \\
   \hspace{15pt} + M_1 + M_1^3 + M_1^3 |t|^4)   
   \hspace{98pt} \textrm{for}\ N=3, \hspace{3pt} (3.54)   \\
   C (\|\Delta \varphi\|_2^2 + \|\nabla \psi_1\|_2^2 + \|\Delta \psi_0\|_2^2 
   + \|\nabla \psi_0\|_2^2   \\
   \hspace{15pt} + M_1 + M_1^3) \exp (C M_1^{1/2} |t|)   
   \hspace{74pt} \textrm{for}\ N=4. \hspace{5pt} (3.55)
   \end{cases}
\end{align*}
The estimate (3.52)-(3.55) are uniform in $n \geq1$ and locally finite in time interval $\subset \mathbb R$.
This implies the global existence of strong solutions $(u_n, v_n, \partial_t v_n)$ as well.
By (3.20) and (3.52)-(3.55), we finally obtain the uniform $H^2 \oplus H^2 \oplus H^1$ bound as
\begin{align*}
   \sup_{n \geq 1} (\|u_n(t)\|_{H^2}^2 + \|v_n(t)\|_{H^2}^2 + \|\partial_t v_n(t)\|_{H^1}^2)
   \leq M_2(t)   \tag{3.56}
\end{align*}
for any $t \in \mathbb R$, where
\begin{align*}
   M_2(t) = 
   \begin{cases}
   C' + C |t|^{4/3} & \textrm{for}\ N=1,   \\
   C' + C t^2 & \textrm{for}\ N=2,   \\
   C' + C t^4 & \textrm{for}\ N=3,   \\
   C' \exp (C |t|) & \textrm{for}\ N=4
   \end{cases}
\end{align*}
with $C$ and $C'$ depending on $\max (\|\varphi\|_{H^1}, \|\psi_0\|_{H^1}, \|\psi_1\|_2)$ and $\max (\|\varphi\|_{H^2}, \|\psi_0\|_{H^2}, \|\psi_1\|_{H^1})$, respectively.

%%%%%% 4. Convergence of the sequence of solutions of the regularized system %%%%%%

\section{Convergence of the sequence of solutions of the regularized system} 

\ \ \ In this section, we prove the convergence of the sequence of solutions $(u_n, v_n, \partial_t v_n)$ of the regularized system (3.1)-(3.3) with values in $H^1 \oplus H^1 \oplus L^2$.
Specifically, we prove that $((u_n, v_n, \partial_t v_n) ; n \geq 1)$ forms a Cauchy sequence in $L^{\infty} (-T, T ; H^1 \oplus H^1 \oplus L^2)$ for any $T>0$.
Below we omit the time variable $t \in [-T,T]$ and we abbreviate $M_2 = M_2(T)$ in (3.56).
For any $m,n \geq 1$ with $m>n$, we start with the equality
\begin{align*}
   &J_m^2 (J_m^2 v_m \cdot J_m^2 u_m) - J_n^2 (J_n^2 v_n \cdot J_n^2 u_n)   \\
   &= J_m^2 ( J_m^2 (v_m - v_n) \cdot J_m^2 u_m 
   + (J_m^2 - J_n^2) v_n \cdot J_m^2 u_m
   + J_n^2 v_n \cdot J_m^2 (u_m - u_n)   \\
   &\hspace{11pt} + J_n^2 v_n \cdot (J_m^2 - J_n^2) u_n )
   + (J_m^2 - J_n^2) (J_n^2 v_n \cdot J_n^2 u_n)   \tag{4.1}
\end{align*}
to compute
\begin{align*}
   &\frac{d}{dt} \|u_m - u_n\|_2^2
   = 2 \textrm{Im} (i \partial_t (u_m - u_n) | u_m - u_n)   \\[-3pt]
   &= 2 \textrm{Im} (- \Delta (u_m - u_n) - J_m^2 (J_m^2 v_m \cdot J_m^2 u_m) 
   - J_n^2 (J_n^2 v_n \cdot J_n^2 u_n) | u_m - u_n)   \\
   &= - 2 \textrm{Im} (J_m^2 (v_m - v_n) \cdot J_m^2 u_m | J_m^2 (u_m - u_n))   \\
   &\hspace{11pt} 
   - 2 \textrm{Im} ((J_m^2 - J_n^2) v_n \cdot J_m^2 u_m | J_m^2 (u_m - u_n))   \\
   &\hspace{11pt} 
   - 2 \textrm{Im} (J_n^2 v_n \cdot (J_m^2 - J_n^2) u_n | J_m^2 (u_m - u_n))   \\
   &\hspace{11pt} 
   - 2 \textrm{Im} ((J_m^2 - J_n^2) (J_n^2 v_n \cdot J_n^2 u_n) | u_m - u_n)   
   =: \textrm{I}_1 + \textrm{I}_2 + \textrm{I}_3 + \textrm{I}_4,   \tag{4.2}
\end{align*}
where we have used
\begin{align*}
   \textrm{Im} (J_m^2 (J_n^2 v_n \cdot J_m^2 (u_m - u_n)) | u_m - u_n)
   = \textrm{Im} (J_n^2 v_n | |J_m^2 (u_m - u_n)|^2) = 0.
\end{align*}
For $(v_n, \partial_t v_n)$, we compute
\begin{align*}
   &\frac{d}{dt} (\|v_m - v_n\|_2^2 + \|\omega^{-1} \partial_t (v_m - v_n)\|_2^2)   \\
   &= 2 (\partial_t(v_m - v_n) | v_m - v_n)
   + (\omega^{-2} \partial_t^2 (v_m - v_n) | \partial_t (v_m - v_n))   \\
   &= 2 (\omega^{-2} (\partial_t^2 - \Delta + 1) (v_m - v_n) | \partial_t (v_m - v_n))   \\
   &= 2 (\omega^{-1} (J_m^2 |J_m^2 u_m|^2 - J_n^2 |J_n^2 u_n|^2) 
   | \omega^{-1} \partial_t (v_m - v_n))   \\
   &= 2 (\omega^{-1} J_m^2 (J_m^2 (u_m - u_n) \cdot \overline{J_m^2 u_m}) 
   | \omega^{-1} \partial_t (v_m - v_n))   \\
   &\hspace{11pt} 
   + 2 (\omega^{-1} J_m^2 ((J_m^2 - J_n^2) u_n \cdot \overline{J_m^2 u_m}) 
   | \omega^{-1} \partial_t (v_m - v_n))   \\
   &\hspace{11pt} 
   + 2 (\omega^{-1} J_m^2 (J_n^2 u_n \cdot \overline{J_m^2 (u_m - u_n)}) 
   | \omega^{-1} \partial_t (v_m - v_n))   \\
   &\hspace{11pt} 
   + 2 (\omega^{-1} J_m^2 (J_n^2 u_n \cdot \overline{(J_m^2 - J_n^2) u_n}) 
   | \omega^{-1} \partial_t (v_m - v_n))   \\
   &\hspace{11pt} 
   + 2 (\omega^{-1} (J_m^2 - J_n^2) |J_n^2 u_n|^2) 
   | \omega^{-1} \partial_t (v_m - v_n))   \\
   &=: \textrm{II}_1 + \textrm{II}_2 + \textrm{II}_3 + \textrm{II}_4 + \textrm{II}_5.   \tag{4.3}
\end{align*}
We estimate terms on the RHS of (4.2) and (4.3) by the H\"{o}lder inequality and Lemma 1-3 to obtain a differential inequality for $\|u_m - u_n\|_2^2 + \|v_m - v_n\|_2^2 + \|\omega^{-1} (v_m - v_n)\|_2^2$.
The term $\textrm{I}_1$ is estimated by
\begin{align*}
   \hspace{50pt} |\textrm{I}_1| 
   &\leq 2 \|v_m - v_n\|_2 \|J_m^2 u_m\|_{\infty} \|u_m - u_n\|_2   \\
   &\leq 
   \begin{cases}
   C M_1^{1/2} (\|u_m - u_n\|_2^2 + \|v_m - v_n\|_2^2) 
   &\hspace{21pt} \textrm{for}\ N=1, \hspace{10pt} (4.4)   \\
   C M_2^{1/2} (\|u_m - u_n\|_2^2 + \|v_m - v_n\|_2^2) 
   &\hspace{21pt} \textrm{for}\ N=2,3, \hspace{1pt} (4.5)  
   \end{cases}
\end{align*}
\vspace{-15pt}
\begin{align*}   
   \hspace{3pt} |\textrm{I}_1| &\leq 2 \|v_m - v_n\|_2 \|J_m^2 u_m\|_4 \|J_m^2 (u_m - u_n)\|_4   \\
   &\leq C \|v_m - v_n\|_2 \|u_m\|_{H^1} \|u_m - u_n\|_{H^1}   \\
   &\leq C M_1^{1/2} (\|\nabla (u_m - u_n)\|_2^2 + \|u_m - u_n\|_2^2 + \|v_m - v_n\|_2^2)
   \hspace{7pt} \textrm{for}\ N=4.   \tag{4.6}
\end{align*}
The term $\textrm{I}_2$ is estimated by
\begin{align*}
   \hspace{50pt} |\textrm{I}_2| 
   &\leq 4 \|(J_m - J_n) v_n\|_2 \|J_m^2 u_m\|_{\infty} \|u_m - u_n\|_2   \\
   &\leq 
   \begin{cases}
   C n^{-1/2} M_1 \|u_m - u_n\|_2
   &\hspace{42pt} \textrm{for}\ N=1, \hspace{10pt} (4.7)   \\
   C n^{-1/2} M_1^{1/2} M_2^{1/2} \|u_m - u_n\|_2
   &\hspace{42pt} \textrm{for}\ N=2,3, \hspace{1pt} (4.8)  
   \end{cases}
\end{align*}
\vspace{-15pt}
\begin{align*}
   \hspace{40pt} |\textrm{I}_2|
   &\leq 2 \|(J_m^2 - J_n^2) v_n\|_4 \|J_m^2 u_m\|_4 \|u_m - u_n\|_2   \\
   &\leq C \|(J_m^2 - J_n^2) v_n\|_{H^1} \|u_m\|_{H^1} \|u_m - u_n\|_2   \\
   &\leq C M_1^{1/2} \|(J_m - J_n) v_n\|_{H^1} \|u_m - u_n\|_2   \\
   &\leq C n^{-1/2} M_1^{1/2} M_2^{1/2} \|u_m - u_n\|_2 
   \hspace{60pt} \textrm{for}\ N=4.   \tag{4.9}
\end{align*}
The term $\textrm{I}_3$ is estimated by
\begin{align*}
   \hspace{42pt} |\textrm{I}_3| 
   &\leq 4 \|J_n^2 v_n\|_{\infty} \|(J_m - J_n) u_n\|_2 \|u_m - u_n\|_2   \\
   &\leq C n^{-1/2} M_1 \|u_m - u_n\|_2 
   \hspace{78pt} \textrm{for}\ N=1,   \tag{4.10}   \\
   |\textrm{I}_3|
   &\leq 2 \|J_n^2 v_n\|_4 \|(J_m^2 - J_n^2) u_n\|_4 \|u_m - u_n\|_2   \\
   &\leq C \|v_n\|_{H^1} \|(J_m - J_n) u_n\|_{H^1} \|u_m - u_n\|_2   \\
   &\leq C n^{-1/2} M_1^{1/2} M_2^{1/2} \|u_m - u_n\|_2 
   \hspace{45pt} \textrm{for}\ N=2,3,4.   \tag{4.11}
\end{align*}
The term $\textrm{I}_4$ is estimated by
\begin{align*}
   \hspace{55pt} |\textrm{I}_4| 
   &\leq 4 \|(J_m - J_n) (J_n^2 v_n \cdot J_n^2 u_n)\|_2 \|u_m - u_n\|_2   \\
   &\leq 4 n^{-1/2} \|\nabla (J_n^2 v_n \cdot J_n^2 u_n)\|_2 \|u_m - u_n\|_2   \\
   &\leq
   \begin{cases}
   C n^{-1/2} M_1 \|u_m - u_n\|_2 
   &\hspace{20pt} \textrm{for}\ N=1, \hspace{22pt} (4.12)   \\
   C n^{-1/2} M_1^{1/2} M_2^{1/2} \|u_m - u_n\|_2 
   &\hspace{20pt} \textrm{for}\ N=2,3,4. \hspace{5pt} (4.13)
   \end{cases}
\end{align*}
The term $\textrm{II}_1$ is estimated by
\begin{align*}
   \hspace{5pt} |\textrm{II}_1| 
   &\leq 2 \|\omega^{-1} (J_m^2 (u_m - u_n) \cdot \overline{J_m^2 u_m} )\|_2 
   \|\omega^{-1} \partial_t (v_m - v_n)\|_2   \\
   &\leq C \|u_m - u_n\|_2 \|J_m^2 u_m\|_{\infty} 
   \|\omega^{-1} \partial_t (v_m - v_n)\|_2   \\
   &\leq
   \begin{cases}
   C M_1^{1/2} (\|u_m - u_n\|_2^2 + \|\omega^{-1} \partial_t (v_m - v_n)\|_2^2) 
   &\hspace{12pt} \textrm{for}\ N=1, \hspace{20pt} (4.14)   \\
   C M_2^{1/2} (\|u_m - u_n\|_2^2 + \|\omega^{-1} \partial_t (v_m - v_n)\|_2^2)
   &\hspace{12pt} \textrm{for}\ N=2,3, \hspace{11pt} (4.15)
   \end{cases}
\end{align*}
\vspace{-14pt}
\begin{align*}
   \hspace{-10pt} |\textrm{II}_1|
   &\leq 2 \|J_m^2 (u_m - u_n)\|_4 \|J_m^2 u_m\|_4 
   \|\omega^{-1} \partial_t (v_m- v_n)\|_2^2   \\
   &\leq C M_2^{1/2} (\|\nabla (u_m - u_n)\|_2^2 + \|u_m - u_n\|_2^2 
   + \|\omega^{-1} \partial_t (v_m - v_n)\|_2^2)  \hspace{30pt}   \\
   &\hspace{230pt} \textrm{for}\ N=4.   \tag{4.16}
\end{align*}
The term $\textrm{II}_2$ is estimated by
\begin{align*}
   \hspace{18pt} |\textrm{II}_2| 
   &\leq 2 \|\omega^{-1} ((J_m^2 - J_n^2) u_n \cdot \overline{J_m^2 u_m} )\|_2 
   \|\omega^{-1} \partial_t (v_m - v_n)\|_2   \\
   &\leq C \|(J_m^2 - J_n^2) u_n\|_2 \|J_m^2 u_m\|_{\infty} 
   \|\omega^{-1} \partial_t (v_m - v_n)\|_2   \\
   &\leq
   \begin{cases}
   C n^{-1/2} M_1 \|\omega^{-1} \partial_t (v_m - v_n)\|_2
   &\hspace{20pt} \textrm{for}\ N=1, \hspace{22pt} (4.17)   \\
   C n^{-1/2} M_1^{1/2} M_2^{1/2} \|\omega^{-1} \partial_t (v_m - v_n)\|_2
   &\hspace{20pt} \textrm{for}\ N=2,3, \hspace{13pt} (4.18)
   \end{cases}
\end{align*}
\vspace{-15pt}
\begin{align*}
   \hspace{-5pt} |\textrm{II}_2|
   &\leq 2 \|(J_m^2 - J_n^2) u_n\|_4 \|J_m^2 u_m\|_4 
   \|\omega^{-1} \partial_t (v_m- v_n)\|_2   \\
   &\leq C n^{-1/2} M_1^{1/2} M_2^{1/2} \|\omega^{-1} \partial_t (v_m - v_n)\|_2
   \hspace{37pt} \textrm{for}\ N=4.   \tag{4.19}
\end{align*}
The term $\textrm{II}_3$ is estimated by
\begin{align*}
   \hspace{5pt} |\textrm{II}_3| 
   &\leq 2 \|\omega^{-1} (J_n^2 u_n \cdot \overline{J_m^2 (u_m - u_n)} )\|_2 
   \|\omega^{-1} \partial_t (v_m - v_n)\|_2   \\
   &\leq C \|J_n^2 u_n\|_{\infty} \|u_m - u_n\|_2 
   \|\omega^{-1} \partial_t (v_m - v_n)\|_2   \\
   &\leq
   \begin{cases}
   C M_1^{1/2} (\|u_m - u_n\|_2^2 + \|\omega^{-1} \partial_t (v_m - v_n)\|_2^2)
   &\hspace{9pt} \textrm{for}\ N=1, \hspace{22pt} (4.20)   \\
   C M_2^{1/2} (\|u_m - u_n\|_2^2 + \|\omega^{-1} \partial_t (v_m - v_n)\|_2^2)
   &\hspace{9pt} \textrm{for}\ N=2,3, \hspace{13pt} (4.21)
   \end{cases}
\end{align*}
\vspace{-15pt}
\begin{align*}
   \hspace{-15pt} |\textrm{II}_3|
   &\leq 2 \|J_n^2 u_n\|_4 \|J_m^2 (u_m - u_n)\|_4 
   \|\omega^{-1} \partial_t (v_m- v_n)\|_2   \\
   &\leq C M_1^{1/2} 
   (\|\nabla (u_m - u_n)\|_2^2 + \|u_m - u_n\|_2^2 
   + \|\omega^{-1} \partial_t (v_m - v_n)\|_2^2)   \\
   &\hspace{225pt} \textrm{for}\ N=4. \hspace{5pt}  \tag{4.22}
\end{align*}
The term $\textrm{II}_4$ is estimated by
\begin{align*}
   \hspace{20pt} |\textrm{II}_4| 
   &\leq 2 \|\omega^{-1} (J_n^2 u_n \cdot \overline{(J_m^2 - J_n^2) u_n} )\|_2 
   \|\omega^{-1} \partial_t (v_m - v_n)\|_2   \\
   &\leq C \|J_n^2 u_n\|_{\infty} \|(J_m^2 - J_n^2) u_n\|_2 
   \|\omega^{-1} \partial_t (v_m - v_n)\|_2   \\
   &\leq
   \begin{cases}
   C n^{-1/2} M_1 \|\omega^{-1} \partial_t (v_m - v_n)\|_2
   &\hspace{19pt} \textrm{for}\ N=1, \hspace{21pt} (4.23)   \\
   C n^{-1/2} M_1^{1/2} M_2^{1/2} \|\omega^{-1} \partial_t (v_m - v_n)\|_2
   &\hspace{19pt} \textrm{for}\ N=2,3, \hspace{12pt} (4.24)
   \end{cases}
\end{align*}
\vspace{-15pt}
\begin{align*}
   \hspace{-2pt} |\textrm{II}_4|
   &\leq 2 \|J_n^2 u_n\|_4 \|(J_m^2 - J_n^2) u_n\|_4 
   \|\omega^{-1} \partial_t (v_m- v_n)\|_2   \\
   &\leq C n^{-1/2} M_1^{1/2} M_2^{1/2} \|\omega^{-1} \partial_t (v_m - v_n)\|_2
   \hspace{37pt} \textrm{for}\ N=4.   \tag{4.25}
\end{align*}
The term $\textrm{II}_5$ is estimated by
\begin{align*}
   \hspace{15pt} |\textrm{II}_5| 
   &\leq 4 \|(J_m - J_n) |J_n u_n|^2\|_2 
   \|\omega^{-1} \partial_t (v_m - v_n)\|_2   \\
   &\leq 4 n^{-1/2} \|\nabla |J_n^2 u_n|^2\|_2  
   \|\omega^{-1} \partial_t (v_m - v_n)\|_2   \\
   &\leq
   \begin{cases}
   C n^{-1/2} M_1 \|\omega^{-1} \partial_t (v_m - v_n)\|_2
   &\hspace{25pt} \textrm{for}\ N=1, \hspace{19pt} (4.26)   \\
   C n^{-1/2} M_1^{1/2} M_2^{1/2} \|\omega^{-1} \partial_t (v_m - v_n)\|_2
   &\hspace{25pt} \textrm{for}\ N=2,3,4. \hspace{2pt} (4.27)
   \end{cases}
\end{align*}
By (4.2)-(4.27), we obtain
\begin{align*}
   &\frac{d}{dt} (\|u_m - u_n\|_2^2 + \|v_m - v_n\|_2^2 
   + \|\omega^{-1} \partial_t (v_m - v_n)\|_2^2)   \\
   &\leq
   \begin{cases}
   C M_1^{1/2} (\|u_m - u_n\|_2^2 + \|v_m - v_n\|_2^2 
   + \|\omega^{-1} \partial_t (v_m - v_n)\|_2^2)   \\
   \ \ \ + C n^{-1/2} M_1 (\|u_m - u_n\|_2 + \|\omega^{-1} \partial_t (v_m - v_n)\|_2) \\
   \hspace{230pt} \textrm{for}\ N=1, \hspace{19pt} (4.28)   
   \\[3pt]
   C M_2^{1/2} (\|u_m - u_n\|_2^2 + \|v_m - v_n\|_2^2 
   + \|\omega^{-1} \partial_t (v_m - v_n)\|_2^2)   \\
   \ \ \ + C n^{-1/2} M_1^{1/2} M_2^{1/2} 
   (\|u_m - u_n\|_2 + \|\omega^{-1} \partial_t (v_m - v_n)\|_2)   \\
   \hspace{230pt} \textrm{for}\ N=2,3, \hspace{9pt} (4.29)   
   \\[3pt]
   C M_2^{1/2} (\|\nabla (u_m - u_n)\|_2^2 + \|u_m - u_n\|_2^2 + \|v_m - v_n\|_2^2 
   + \|\omega^{-1} \partial_t (v_m - v_n)\|_2^2)   \\
   \ \ \ + C n^{-1/2} M_1^{1/2} M_2^{1/2} (\|u_m - u_n\|_2 
   + \|\omega^{-1} \partial_t (v_m - v_n)\|_2) \\
   \hspace{230pt} \textrm{for}\ N=4. \hspace{19pt} (4.30)   
   \end{cases}
\end{align*}
The differential inequalities (4.28) and (4.29) are closed and integrated by the Gronwall lemma as
\begin{align*}
   &\|(u_m - u_n)(t)\|_2^2 + \|(v_m - v_n)(t)\|_2^2 
   + \|\omega^{-1} \partial_t (v_m - v_n)(t)\|_2^2   \\
   &\leq 
   \begin{cases}
   \exp (C M_1^{1/2} t) (\|(J_m^2 - J_n^2) \varphi\|_2^2 
   + \|(J_m^2 - J_n^2) \psi_0\|_2^2 
   + \|\omega^{-1} (J_m^2 - J_n^2) \psi_1\|_2^2)   \\
   \ \ \ + C n^{-1/2} M_1^{3/2} t
   \hspace{151pt} \textrm{for}\ N=1, \hspace{18pt} (4.31)   
   \\[3pt]
   \exp (C M_1^{1/2} t) (\|(J_m^2 - J_n^2) \varphi\|_2^2 
   + \|(J_m^2 - J_n^2) \psi_0\|_2^2 
   + \|\omega^{-1} (J_m^2 - J_n^2) \psi_1\|_2^2)   \\
   \ \ \ + C n^{-1/2} M_1^{1/2} M_2 t
   \hspace{137pt} \textrm{for}\ N=2,3. \hspace{10pt} (4.32)   
   \end{cases}
\end{align*}
It follows from (4.31) and (4.32) that $((u_n, v_n, \partial_t v_n) ; n \geq 1)$ is a Cauchy sequence in $L^{\infty}(-T,T ; L^2 \oplus L^2 \oplus H^{-1})$ for any $T>0$ if $N \leq 3$.
By the boundedness in $L^{\infty}(-T,T; H^2 \oplus H^2 \oplus H^1)$, this implies that $((u_n, v_n, \partial_t v_n) ; n \geq 1)$ is also a Cauchy sequence in $L^{\infty}(-T,T ; H^1 \oplus H^1 \oplus L^2)$ for any $T>0$.

Therefore, it remains to prove the convergence of $((u_n, v_n, \partial_t v_n) ; n \geq 1)$ in $L^{\infty}(-T,T ; H^1 \oplus H^1 \oplus L^2)$ for any $T>0$ in the case $N=4$.
For that purpose we compute
\begin{align*}
   &\frac{d}{dt} (\|\nabla (u_m - u_n)\|_2^2 
   + \frac{1}{2} (\|\nabla (v_m - v_n)\|_2^2 + \|v_m - v_n\|_2^2 
   + \|\partial_t (v_m - v_n)\|_2^2) )   \\
   &= 2 \textrm{Re} (\partial_t \nabla (u_m - u_n) | \nabla (u_m - u_n))
   + (\partial_t \nabla (v_m - v_n) | \nabla (v_m - v_n))   \\
   &\hspace{11pt} + (\partial_t (v_m - v_n) | v_m - v_n)
   + (\partial_t^2 (v_m - v_n) | \partial_t (v_m - v_n))   \\
   &= - 2 \textrm{Re} (\partial_t (u_m - u_n) | \Delta (u_m - u_n))
   + ((\partial_t^2 - \Delta + 1) (v_m - v_n) | \partial_t (v_m - v_n))   \\
   &= - 2 \textrm{Re} (\partial_t (u_m - u_n) | - i \partial_t (u_m - u_n) 
   - J_m^2 (J_m^2 v_m \cdot J_m^2 u_m^2) + J_n^2 (J_n^2 v_n \cdot J_n^2 u_n))   \\
   &\hspace{11pt} + (J_m^2 |J_m^2 u_m|^2 - J_n^2 |J_n^2 u_n|^2 | \partial_t (v_m - v_n)) \\
   &= 2 \textrm{Re} (J_m^2 (J_m^2 v_m \cdot J_m^2 u_m^2) 
   - J_n^2 (J_n^2 v_n \cdot J_n^2 u_n) | \partial_t (u_m - u_n))   \\
   &\hspace{11pt} + (J_m^2 |J_m^2 u_m|^2 - J_n^2 |J_n^2 u_n|^2 | \partial_t (v_m - v_n)) \\
   &= 2 \textrm{Re} (J_m^2 ( J_m^2 (v_m - v_n) \cdot J_m^2 (u_m - u_n)
   + J_m^2 (v_m - v_n) \cdot J_m^2 u_n   \\ 
   &\hspace{50pt} + J_n^2 v_n \cdot J_m^2 (u_m - u_n)   
   + (J_m^2 - J_n^2) v_n \cdot J_m^2 u_m   \\
   &\hspace{50pt} + J_n^2 v_n \cdot (J_m^2 - J_n^2) u_n )
   + (J_m^2 - J_n^2) (J_n^2 v_n \cdot J_n^2 u_n) | \partial_t (u_m - u_n))   \\
   &\hspace{11pt} + \textrm{Re} (J_m^2 (|J_m^2 (u_m - u_n)|^2 
   + J_m^2 (u_m - u_n) \cdot \overline{J_m^2 u_n}  
   + (J_m^2 - J_n^2) u_n \cdot \overline{J_m^2 u_m}   \\
   &\hspace{50pt} + J_n^2 u_n \cdot \overline{J_m^2 (u_m - u_n)} + J_n^2 u_n \cdot \overline{(J_m^2 - J_n^2) u_n})   \\
   &\hspace{50pt} + (J_m^2 - J_n^2) |J_n u_n|^2 | \partial_t (v_m - v_n) )   \\
   &= \frac{d}{dt} (|J_m^2 (u_m - u_n)|^2 | J_m^2 (v_m - v_n))   \\
   &\hspace{11pt} + 2 \textrm{Re} (J_m^2 (J_m^2 (v_m - v_n) \cdot J_m^2 u_n) 
   | \partial_t (u_m - u_n))   \\
   &\hspace{11pt} + 2 \textrm{Re} (J_m^2 (J_m^2 v_n \cdot J_m^2 (u_m - u_n)) 
   | \partial_t (u_m - u_n))   \\
   &\hspace{11pt} + 2 \textrm{Re} (J_m^2 ((J_m^2 - J_n^2) v_n \cdot J_m^2 u_m) 
   | \partial_t (u_m - u_n))   \\
   &\hspace{11pt} + 2 \textrm{Re} (J_m^2 (J_m^2 v_n \cdot (J_m^2 - J_n^2) u_n) 
   | \partial_t (u_m - u_n))   \\
   &\hspace{11pt} + 2 \textrm{Re} ((J_m^2 - J_n^2) (J_n^2 v_n \cdot J_n^2 u_n) 
   | \partial_t (u_m - u_n))   \\
   &\hspace{11pt} + \textrm{Re} (J_m^2 (J_m^2 (u_m - u_n) \cdot \overline{J_m u_m}) 
   | \partial_t (v_m - v_n))   \\
   &\hspace{11pt} + \textrm{Re} (J_m^2 ((J_m^2 - J_n^2) u_n \cdot \overline{J_m^2 u_n}) 
   | \partial_t (v_m - v_n))   \\
   &\hspace{11pt} + \textrm{Re} (J_m^2 ( J_n^2 u_n \cdot \overline{J_m^2 (u_m - u_n)}) 
   | \partial_t (v_m - v_n))   \\
   &\hspace{11pt} + \textrm{Re} (J_m^2 (J_n^2 u_n \cdot \overline{(J_m^2 - J_n^2) u_n}) 
   | \partial_t (v_m - v_n))   \\
   &\hspace{11pt} + \textrm{Re} ( (J_m^2 - J_n^2) |J_n^2 u_n|^2 
   | \partial_t (v_m - v_n))   
\end{align*}     
\begin{align*}
   &= \frac{d}{dt} (|J_m^2 (u_m - u_n)|^2 | J_m (v_m - v_n))   \\
   &\hspace{11pt} + \frac{d}{dt} (2 \textrm{Re} (J_m^2 (v_m - v_n) \cdot J_m^2 u_n 
   | J_m^2 (u_m - u_n)))   \\
   &\hspace{26pt} - 2 \textrm{Re} (J_m^2 \partial_t (v_m - v_n) \cdot J_m^2 u_n 
   | J_m^2 (u_m - u_n))   \\
   &\hspace{26pt} - 2 \textrm{Re} (J_m^2 (v_m - v_n) \cdot J_m^2 \partial_t u_n
   | J_m^2 (u_m - u_n))   \\
   &\hspace{11pt} + \frac{d}{dt} (J_m^2 v_n  | |J_m^2 (u_m - u_n)|^2)   
   - (J_m^2 \partial_t v_n  | |J_m^2 (u_m - u_n)|^2)   \\
   &\hspace{11pt} + \frac{d}{dt} (2 \textrm{Re} ((J_m^2 - J_n^2) v_n \cdot J_m^2 u_m 
   | J_m^2 (u_m - u_n)))   \\
   &\hspace{26pt} - 2 \textrm{Re} ((J_m^2 - J_n^2) \partial_t v_n \cdot J_m^2 u_m 
   | J_m^2 (u_m - u_n))   \\
   &\hspace{26pt} - 2 \textrm{Re} ((J_m^2 - J_n^2) v_n \cdot J_m^2 \partial_t u_m
   | J_m^2 (u_m - u_n))   \\
   &\hspace{11pt} + \frac{d}{dt} (2 \textrm{Re} (J_m^2 v_n \cdot (J_m^2 - J_n^2) u_n 
   | J_m^2 (u_m - u_n)))   \\
   &\hspace{26pt} - 2 \textrm{Re} (J_m^2 \partial_t v_n \cdot (J_m^2 - J_n^2) u_n 
   | J_m^2 (u_m - u_n))   \\
   &\hspace{26pt} - 2 \textrm{Re} (J_m^2 v_n \cdot (J_m^2 - J_n^2) \partial_t u_n
   | J_m^2 (u_m - u_n))   \\
   &\hspace{11pt} + \frac{d}{dt} (2 \textrm{Re} ((J_m^2 - J_n^2) (J_n^2 v_n \cdot J_n^2 u_n) 
   | u_m - u_n))   \\
   &\hspace{26pt} - 2 \textrm{Re} ((J_m^2 - J_n^2) (J_n^2 \partial_t v_n \cdot J_n^2 u_n) 
   | u_m - u_n)   \\
   &\hspace{26pt} - 2 \textrm{Re} ((J_m^2 - J_n^2) (J_n^2 v_n \cdot J_n^2 \partial_t u_n) 
   | u_m - u_n)   \\
   &\hspace{11pt} + \textrm{Re} (J_m^2 (J_m^2 (u_m - u_n) \cdot \overline{J_m u_m}) 
   | \partial_t (v_m - v_n))   \\
   &\hspace{11pt} + \textrm{Re} (J_m^2 ((J_m^2 - J_n^2) u_n \cdot \overline{J_m^2 u_n}) 
   | \partial_t (v_m - v_n))   \\
   &\hspace{11pt} + \textrm{Re} (J_m^2 ( J_n^2 u_n \cdot \overline{J_m^2 (u_m - u_n)}) 
   | \partial_t (v_m - v_n))   \\
   &\hspace{11pt} + \textrm{Re} (J_m^2 (J_n^2 u_n \cdot \overline{(J_m^2 - J_n^2) u_n}) 
   | \partial_t (v_m - v_n))   \\
   &\hspace{11pt} + \textrm{Re} ( (J_m^2 - J_n^2) |J_n^2 u_n|^2 
   | \partial_t (v_m - v_n)).   \tag{4.33}  
\end{align*}
For any $m, n$ with $m>n \geq 1$, we introduce $E_{m,n}(t)$ by
\begin{align*}
   E_{m,n}(t) &:= \|\nabla (u_m - u_n)\|_2^2 
   + \frac{1}{2} (\|\nabla (v_m - v_n)\|_2^2 + \|v_m - v_n\|_2^2 
   + \|\partial_t (v_m - v_n)\|_2^2)   \\
   &\hspace{13pt} - (|J_m^2 (u_m - u_n)|^2 | J_m (v_m - v_n))   \\
   &\hspace{13pt} - 2 \textrm{Re} (J_m^2 (v_m - v_n) \cdot J_m^2 u_n 
   | J_m^2 (u_m - u_n))   \\
   &\hspace{13pt} - (J_m^2 v_n  | |J_m^2 (u_m - u_n)|^2)   \\
   &\hspace{13pt} - 2 \textrm{Re} ((J_m^2 - J_n^2) v_n \cdot J_m^2 u_m 
   | J_m^2 (u_m - u_n))   \\
   &\hspace{13pt} - 2 \textrm{Re} (J_m^2 v_n \cdot (J_m^2 - J_n^2) u_n 
   | J_m^2 (u_m - u_n))   \\
   &\hspace{13pt} - 2 \textrm{Re} ((J_m^2 - J_n^2) (J_n^2 v_n \cdot J_n^2 u_n) 
   | u_m - u_n).   \tag{4.34}
\end{align*}
Then the equality (4.33) is rewritten as
\begin{align*}
   E_{m,n}'(t)
   &= - 2 \textrm{Re} (J_m^2 \partial_t (v_m - v_n) \cdot J_m^2 u_n 
   | J_m^2 (u_m - u_n))   \\ 
   &\hspace{11pt} - 2 \textrm{Re} (J_m^2 (v_m - v_n) \cdot J_m^2 \partial_t u_n
   | J_m^2 (u_m - u_n))   \\
   &\hspace{11pt} - (J_m^2 \partial_t v_n  | |J_m^2 (u_m - u_n)|^2)   \\
   &\hspace{11pt} - 2 \textrm{Re} ((J_m^2 - J_n^2) \partial_t v_n \cdot J_m^2 u_n 
   | J_m^2 (u_m - u_n))   \\
   &\hspace{11pt} - 2 \textrm{Re} ((J_m^2 - J_n^2) v_n \cdot J_m^2 \partial_t u_n
   | J_m^2 (u_m - u_n))   \\
   &\hspace{11pt} - 2 \textrm{Re} (J_m^2 \partial_t v_n \cdot (J_m^2 - J_n^2) u_n 
   | J_m^2 (u_m - u_n))   \\
   &\hspace{11pt} - 2 \textrm{Re} (J_m^2 v_n \cdot (J_m^2 - J_n^2) \partial_t u_n
   | J_m^2 (u_m - u_n))   \\
   &\hspace{11pt} - 2 \textrm{Re} ((J_m^2 - J_n^2) (J_n^2 \partial_t v_n \cdot J_n^2 u_n) 
   | J_m^2 (u_m - u_n))   \\
   &\hspace{11pt} - 2 \textrm{Re} ((J_m^2 - J_n^2) (J_n^2 v_n \cdot J_n^2 \partial_t u_n) 
   | J_m^2 (u_m - u_n))   \\
   &\hspace{11pt} + \textrm{Re} (J_m^2 (J_m^2 (u_m - u_n) \cdot \overline{J_m u_m}) 
   | \partial_t (v_m - v_n))   \\
   &\hspace{11pt} + \textrm{Re} (J_m^2 ((J_m^2 - J_n^2) u_n \cdot \overline{J_m^2 u_n}) 
   | \partial_t (v_m - v_n))   \\
   &\hspace{11pt} + \textrm{Re} (J_m^2 ( J_n^2 u_n \cdot \overline{J_m^2 (u_m - u_n)}) 
   | \partial_t (v_m - v_n))   \\
   &\hspace{11pt} + \textrm{Re} (J_m^2 (J_n^2 u_n \cdot \overline{(J_m^2 - J_n^2) u_n}) 
   | \partial_t (v_m - v_n))   \\
   &\hspace{11pt} + \textrm{Re} ( (J_m^2 - J_n^2) |J_n^2 u_n|^2 
   | \partial_t (v_m - v_n))   \\
   &=: \textrm{III}_1 + \cdots + \textrm{III}_{14}.   \tag{4.35}
\end{align*}
We estimate indefinite terms on the RHS of (4.34) and (4.35) by the H\"{o}lder inequality and Lemmas 1-3 to obtain a differential inequality for $E_{m,n}$.

The first and the third indefinite terms on the RHS of (4.34) is estimated by
\begin{align*}
   &|- (|J_m^2 (u_m - u_n)|^2 | J_m (v_m - v_n)) - (J_m^2 v_n  | |J_m^2 (u_m - u_n)|^2)| \\
   &= |(J_m^2 v_m | |J_m^2 (u_m - u_n)|^2)|   \\
   &\leq \|J_m^2 v_m\|_4 \|J_m^2 (u_m - u_n)\|_{8/3}^2   \\
   &\leq C \|v_m\|_{H^1} \|u_m - u_n\|_{H^1} \|u_m - u_n\|_2   \\
   &\leq C M_2^{1/2} \|u_m - u_n\|_2 \|\nabla (u_m - u_n)\|_2 
   + C M_2^{1/2} \|u_m - u_n\|_2^2   \\
   &\leq \varepsilon \|\nabla (u_m - u_n)\|_2^2 
   + C (\varepsilon^{-1} M_2 + M_2^{1/2}) \|u_m - u_n\|_2^2,   \tag{4.36}
\end{align*}
where $\varepsilon > 0$ sufficiently small.
The second indefinite term on the RHS of (4.34) is estimated by
\begin{align*}
   &|2 \textrm{Re} (J_m^2 (v_m - v_n) \cdot J_m^2 u_n | J_m^2 (u_m - u_n))|   \\
   &\leq 2 C_{4,4}^2 \|u_n\|_2 \|v_m - v_n\|_{H^1} \|u_m - u_n\|_{H^1}   \\
   &\leq 2 C_{4,4}^2 \|\varphi\|_2 (\|\nabla (v_m - v_n)\|_2^2 + \|v_m - v_n\|_2^2)^{1/2}
   (\|\nabla (u_m - u_n)\|_2^2 + \|u_m - u_n\|_2^2)^{1/2}   \\
   &\leq \sqrt{2} C_{4,4}^2 \|\varphi\|_2 
   ( \|\nabla (u_m - u_n)\|_2^2 + \|u_m - u_n\|_2^2   \\
   &\hspace{70pt} + \frac{1}{2} (\|\nabla (v_m - v_n)\|_2^2 + \|v_m - v_n\|_2^2) ). \tag{4.37}
\end{align*}
The fourth indefinite term on the RHS of (4.34) is estimated by
\begin{align*}
   &|2 \textrm{Re} ((J_m^2 - J_n^2) v_n \cdot J_m^2 u_m | J_m^2 (u_m - u_n))|   \\
   &\leq 2 \|(J_m^2 - J_n^2) v_n\|_4 \|J_m^2 u_m\|_4 \|u_m - u_n\|_2   \\
   &\leq C n^{-1/2} M_2 \|u_m - u_n\|_2.   \tag{4.38}
\end{align*}
The fifth indefinite term on the RHS of (4.34) is estimated by
\begin{align*}
   &|2 \textrm{Re} (J_m^2 v_n \cdot (J_m^2 - J_n^2) u_n | J_m^2 (u_m - u_n))|   \\
   &\leq 2 \|J_n^2 v_n\|_4 \|(J_m^2 - J_n^2) u_n\|_4 \|u_m - u_n\|_2   \\
   &\leq C n^{-1/2} M_2 \|u_m - u_n\|_2.   \tag{4.39}
\end{align*}
The last indefinite term on the RHS of (4.34) is estimated by
\begin{align*}
   &|2 \textrm{Re} ((J_m^2 - J_n^2) (J_n^2 v_n \cdot J_n^2 u_n) | u_m - u_n)|   \\
   &\leq 2 \|(J_m^2 - J_n^2) (J_n^2 v_n \cdot J_n^2 u_n)\|_2 \|u_m - u_n\|_2   \\
   &\leq 2 n^{-1/2} \|\nabla (J_n^2 v_n \cdot J_n^2 u_n)\|_2 \|u_m - u_n\|_2   \\
   &\leq C n^{-1/2} M_2 \|u_m - u_n\|_2.   \tag{4.40}
\end{align*}
By (4.34) and (4.36)-(4.40), we obtain
\begin{align*}
   &(1 - \varepsilon - \sqrt{2} C_{4,4}^2 \|\varphi\|_2)
   (\|\nabla (u_m - u_n)\|_2^2 
   + \frac{1}{2} (\|\nabla (v_m - v_n)\|_2^2 + \|v_m - v_n\|_2^2))   \\
   &\hspace{11pt} + \frac{1}{2} \|\partial_t (v_m - v_n)\|_2^2   \\
   &\leq E_{m,n}(t) + C (\varepsilon^{-1} M_2 + M_2^{1/2} + \|\varphi\|_2) \|u_m - u_n\|_2^2
   + C n^{-1/2} M_2 \|u_m - u_n\|_2   \\
   &\leq C (M_2^{1/2} + \|\varphi\|_2 + \varepsilon)
   (\|\nabla (u_m - u_n)\|_2^2   \\
   &\hspace{11pt} + \frac{1}{2} (\|\nabla (v_m - v_n)\|_2^2 
   + \|v_m - v_n\|_2^2 + \|\partial_t (v_m - v_n)\|_2^2))   \\
   &\hspace{11pt} + C (\varepsilon^{-1} M_2 + M_2^{1/2} + \|\varphi\|_2) \|u_m - u_n\|_2^2
   + C n^{-1/2} M_2 \|u_m - u_n\|_2.   \tag{4.41}
\end{align*}
If we assume $\sqrt{2} C_{4,4}^2 \|\varphi\|_2 < 1$, we choose $\varepsilon >0$ with $1 - \varepsilon - \sqrt{2} C_{4,4}^2 \|\varphi\|_2 >0$ and the first factor appearing in (4.41) is strictly positive.

We continue to estimate indefinite terms.
The term $\textrm{III}_1$ is estimated by
\begin{align*}
   &|2 \textrm{Re} (J_m^2 \partial_t (v_m - v_n) \cdot J_m^2 u_n| J_m^2 (u_m - u_n))| \\
   &\leq 2 \|J_m^2 \partial_t (v_m - v_n)\|_2 \|J_m^2 u_n\|_4 \|u_m - u_n\|_4   \\
   &\leq C M_2^{1/2} \|\partial_t (v_m - v_n)\|_2 \|u_m - u_n\|_{H^1}.   \tag{4.42}
\end{align*}
The term $\textrm{III}_2$ is estimated by
\begin{align*}
   &|2 \textrm{Re} (J_m^2 (v_m - v_n) \cdot J_m^2 \partial_t u_n| J_m^2 (u_m - u_n))| \\
   &\leq 2 \|J_m^2 (v_m - v_n)\|_4 \|J_m^2 \partial_t u_n\|_4 \|u_m - u_n\|_2   \\
   &\leq C M_2^{1/2} \|v_m - v_n\|_{H^1} \|u_m - u_n\|_2.   \tag{4.43}
\end{align*}
The term $\textrm{III}_3$ is estimated by
\begin{align*}
   &|(J_m^2 \partial_t v_n  | |J_m^2 (u_m - u_n)|^2)| \\
   &\leq \|\partial_t v_m\|_4 \|J_m^2 (u_m - u_n)\|_{8/3}^2   \\
   &\leq C M_2^{1/2} \|u_m - u_n\|_{H^1} \|u_m - u_n\|_2.   \tag{4.44}
\end{align*}
The term $\textrm{III}_4$ is estimated by
\begin{align*}
   &|2 \textrm{Re} ((J_m^2 - J_n^2) \partial_t v_n \cdot J_m^2 u_m
   | J_m^2 (u_m - u_n))| \\
   &\leq 2 \|(J_m^2 - J_n^2) \partial_t v_n\|_4 \|J_m^2 u_m\|_4 \|u_m - u_n\|_2   \\
   &\leq C n^{-1/2} M_2 \|u_m - u_n\|_{H^1}.   \tag{4.45}
\end{align*}
The term $\textrm{III}_5$ is estimated by
\begin{align*}
   &|2 \textrm{Re} ((J_m^2 - J_n^2) v_n \cdot J_m^2 \partial_t u_m
   | J_m^2 (u_m - u_n))| \\
   &\leq 2 \|(J_m^2 - J_n^2) v_n\|_4 \|J_m^2 \partial_t u_m\|_2 \|J_m^2 (u_m - u_n)\|_4  \\
   &\leq C n^{-1/2} M_2 \|u_m - u_n\|_{H^1}.   \tag{4.46}
\end{align*}
The term $\textrm{III}_6$ is estimated by
\begin{align*}
   &|2 \textrm{Re} (J_m^2 \partial_t v_n \cdot (J_m^2 - J_n^2) u_n 
   | J_m^2 (u_m - u_n))| \\
   &\leq 2 \|J_m^2 \partial_t v_n\|_4 \|(J_m^2 - J_n^2) u_n\|_4 \|u_m - u_n\|_2  \\
   &\leq C n^{-1/2} M_2 \|u_m - u_n\|_2.   \tag{4.47}
\end{align*}
The term $\textrm{III}_7$ is estimated by
\begin{align*}
   &|2 \textrm{Re} (J_m^2 v_n \cdot (J_m^2 - J_n^2) \partial_t u_n
   | J_m^2 (u_m - u_n))| \\
   &= |2 \textrm{Re} ( \partial_t u_n
   | (J_m^2 - J_n^2) ( J_m^2 v_n \cdot J_m^2 (u_m - u_n)))| \\
   &\leq 2 \|\partial_t u_n\|_2 \|(J_m^2 - J_n^2) ( J_m^2 v_n \cdot J_m^2 (u_m - u_n))\|_2   \\
   &\leq C n^{-1/2} \|\partial_t u_n\|_2 \|\nabla ( J_m^2 v_n \cdot J_m^2 (u_m - u_n))\|_2   \\
   &\leq C n^{-1/2} \|\partial_t u_n\|_2 ( \|J_m^2 \nabla v_n \cdot J_m^2 (u_m - u_n))\|_2 + \|J_m^2 v_n \cdot J_m^2 \nabla (u_m - u_n))\|_2 )   \\
   &\leq C n^{-1/2} M_2^{3/2}.   \tag{4.48}
\end{align*}
The term $\textrm{III}_8$ is estimated by
\begin{align*}
   &|2 \textrm{Re} ((J_m^2 - J_n^2) (J_n^2 \partial_t v_n \cdot J_n^2 u_n) 
   | u_m - u_n)| \\
   &= |2 \textrm{Re} (J_n^2 \partial_t v_n \cdot J_n^2 u_n 
   | (J_m^2 - J_n^2) (u_m - u_n))|   \\
   &\leq 2 \|J_n^2 \partial_t v_n \cdot J_n^2 u_n\|_2 \|(J_m^2 - J_n^2) (u_m - u_n)\|_2 \\
   &\leq 2 \|J_n \partial_t v_n\|_4 \|J_n u_n\|_4 n^{-1/2} \|\nabla (u_m - u_n)\|_2  \\
   &\leq C n^{-1/2} M_2 \|\nabla (u_m - u_n)\|_2.   \tag{4.49}
\end{align*}
The term $\textrm{III}_9$ is estimated by
\begin{align*}
   &|2 \textrm{Re} ((J_m^2 - J_n^2) (J_n^2 v_n \cdot J_n^2 \partial_t u_n) 
   | u_m - u_n)| \\
   &= |2 \textrm{Re} (J_n^2 v_n \cdot J_n^2 \partial_t u_n 
   | (J_m^2 - J_n^2) (u_m - u_n))|   \\
   &\leq 2 \|J_n^2 v_n\|_4 \|\partial_t u_n\|_2 \|(J_m^2 - J_n^2) (u_m - u_n)\|_4 \\
   &\leq C M_2 \|(J_m^2 - J_n^2) (u_m - u_n)\|_{H^1} \\
   &\leq C n^{-1/2} M_2 \|\nabla (u_m - u_n)\|_{H^1}  \\
   &\leq C n^{-1/2} M_2^{3/2}.   \tag{4.50}
\end{align*}
The term $\textrm{III}_{10}$ is estimated by
\begin{align*}
   &|\textrm{Re} (J_m^2 (u_m - u_n) \cdot \overline{J_m^2 u_m} 
   | J_m^2 \partial_t (v_m - v_n))| \\
   &\leq \|J_m^2 (u_m - u_n)\|_4 \|J_m^2 u_m\|_4 \|\partial_t (v_m - v_n)\|_2 \\
   &\leq C M_1^{1/2} \|u_m - u_n\|_{H^1} \|\partial_t (v_m - v_n)\|_2.   \tag{4.51}
\end{align*}
The term $\textrm{III}_{11}$ is estimated by
\begin{align*}
   &|\textrm{Re} ((J_m^2 - J_n^2) u_n \cdot \overline{J_m^2 u_n} 
   | J_m^2 \partial_t (v_m - v_n))| \\
   &\leq \|J_m^2 (u_m - u_n)\|_4 \|J_m^2 u_n\|_4 \|\partial_t (v_m - v_n)\|_2 \\
   &\leq C M_1^{1/2} \|u_m - u_n\|_{H^1} \|\partial_t (v_m - v_n)\|_2.   \tag{4.52}
\end{align*}
The term $\textrm{III}_{12}$ is estimated by
\begin{align*}
   &|\textrm{Re} (J_n^2 u_n \cdot \overline{J_m^2 (u_m - u_n)} | J_m^2 \partial_t (v_m - v_n))|   \\
   &\leq \|J_n^2 u_n\|_4 \|u_m - u_n\|_4 \|\partial_t (v_m - v_n)\|_2   \\
   &\leq C M_1^{1/2} \|u_m - u_n\|_{H^1} \|\partial_t (v_m - v_n)\|_2.   \tag{4.53}
\end{align*}
The term $\textrm{III}_{13}$ is estimated by
\begin{align*}
   &|\textrm{Re} (J_n^2 u_n \cdot \overline{(J_m^2 - J_n^2) u_n} 
   | J_m^2 \partial_t (v_m - v_n))| \\
   &\leq \|J_n^2 u_n\|_4 \|(J_m^2 - J_n^2) u_n\|_4 \|\partial_t (v_m - v_n)\|_2 \\
   &\leq C n^{-1/2} M_2 \|\partial_t (v_m - v_n)\|_2.   \tag{4.54}
\end{align*}
The term $\textrm{III}_{14}$ is estimated by
\begin{align*}
   &|\textrm{Re} ((J_m^2 - J_n^2) |J_n^2 u_n|^2 | J_m^2 \partial_t (v_m - v_n))| \\
   &\leq \|(J_m^2 - J_n^2) |J_n^2 u_n|^2\|_2 \|\partial_t (v_m - v_n)\|_2   \\
   &\leq n^{-1/2} \|\nabla |J_n^2 u_n|^2\|_2 \|\partial_t (v_m - v_n)\|_2   \\
   &\leq 2 n^{-1/2} \|\nabla J_n^2 u_n\|_4 \|J_n^2 u_n\|_4 \|\partial_t (v_m - v_n)\|_2 \\
   &\leq C n^{-1/2} M_2 \|\partial_t (v_m - v_n)\|_2.   \tag{4.55}
\end{align*}
By (4.35), (4.42)-(4.55), and (4.41), we obtain
\begin{align*}
   |E_{m,n}'(t)| 
   &\leq C M_2^{1/2} (\|\nabla (u_m - u_n)\|_2^2   \\
   &\hspace{11pt} + \frac{1}{2} (\|\nabla (v_m - v_n)\|_2^2 + \|v_m - v_n\|_2^2 
   + \|\partial_t (v_m - v_n)\|_2^2) )   \\
   &\hspace{11pt} + C M_2^{1/2} \|u_m - u_n\|_2^2 + C n^{-1/2} M_2^{3/2}   \\
   &\hspace{11pt} + C n^{-1/2} M_2 ( \|u_m - u_n\|_{H^1}^2 
   + \|\partial_t (v_m - v_n)\|_2^2 )^{1/2}   \\
   &\leq C E_{m,n}(t) + C n^{-1/2} ( E_{m,n}(t) + \|u_m - u_n\|_2^2 + n^{-1} )^{1/2}   \\
   &\hspace{11pt} + C \|u_m - u_n\|_2^2 + C n^{-1/2},   \tag{4.56}
\end{align*}
where we have omitted the dependence on $M_1, M_2$, and $\delta$ on the constants for simplicity.
We now introduce $F_{m,n}(t)$ as
\begin{align*}
   F_{m,n}(t) := E_{m,n}(t) + \|u_m - u_n\|_2^2 + \|v_m - v_n\|_2^2 
   + \|\omega^{-1} \partial_t (v_m - v_n)\|_2^2 + n^{-1}.
\end{align*}
Then it follows from (4.30) and (4.54) that
\begin{align*}
   |F_{m,n}'(t)| \leq C F_{m,n}(t) + C n^{-1/2}.   \tag{4.57}
\end{align*}
By the Gronwall lemma, we have
\begin{align*}
   F_{m,n}(t) \leq e^{C |t|} ( F_{m,n}(0) + C n^{-1/2} |t| ).
\end{align*}
This implies
\begin{align*}
   \sup_{|t| \leq T} (\|u_m(t) - u_n(t)\|_{H^1} + \|v_m(t) - v_n(t)\|_{H^1} 
   + \|\partial_t (v_m(t) - v_n(t))\|_2) \to 0
\end{align*}
as $m > n \to \infty$, as required.

%%%%%% 5. Proofs of Theorems 1 and 2 %%%%%%%

\section{Proofs of Theorems 1 and 2}

\textbf{Proof of Theorem 1.}   
By (4.56), there exists $(u, v, w) \in C(\mathbb R ; H_0^1 \oplus H_0^1 \oplus L^2)$ such that
\begin{align*}
   \lim_{n \to \infty} \sup_{|t| \leq T} (\|u_n(t) - u(t)\|_{H^1} + \|v_n(t) - v(t)\|_{H^1} 
   + \|\partial_t v_n(t) - w(t)\|_2) = 0   \tag{5.1}
\end{align*}
for any $T>0$.
Taking the limit of both sides of the equality
\begin{align*}
   v_n(t) = J_n^2 \psi_0 + \int_0^t \partial_t v_n(s) ds,
\end{align*}
we have
\begin{align*}
   v(t) = \psi_0 + \int_0^t w(s) ds
\end{align*}
so that $v \in C^1 (\mathbb R ; L^2)$ and $\partial_t v = w$.

Part (2) of Theorem 1 follows from (3.6) and (3.7) since
\begin{align*}
   &\sup_{|t| \leq T} |(|J_n^2 u_n(t)|^2 | J_n^2 v_n(t)) - (|u(t)|^2 | v(t))|   \\
   &\leq \sup_{|t| \leq T} |(|J_n^2 u_n(t)|^2 - |u(t)|^2 | J_n^2 v_n(t))|
   + \sup_{|t| \leq T} |(|u(t)|^2 | J_n^2 v_n(t) - v(t))|   \\
   &\leq \sup_{|t| \leq T} (\|J_n^2 u_n(t)\|_4 + \|u(t)\|_4)
   \|J_n^2 u_n(t) - u(t)\|_4 \|J_n^2 v_n(t)\|_2   \\
   &\hspace{11pt} + \sup_{|t| \leq T} \|u(t)\|_4^2 \|J_n^2 v_n(t) - v(t)\|_2   \\
   &\leq C M_1 \sup_{|t| \leq T} (\|J_n^2 (u_n(t) - u(t))\|_{H^1}
   + \|(J_n^2 - I) u(t)\|_{H^1})    \\
   &\hspace{11pt} + C M_1 \sup_{|t| \leq T} (\|J_n^2 (v_n(t) - v(t))\|_2
   + \|(J_n^2 - I) v(t)\|_2)   \\
   &\leq C M_1 \sup_{|t| \leq T} (\|u_n(t) - u(t)\|_{H^1}
   + \|(J_n^2 - I) u(t)\|_{H^1})    \\
   &\hspace{11pt} + C M_1 \sup_{|t| \leq T} (\|v_n(t) - v(t)\|_2
   + \|(J_n^2 - I) v(t)\|_2) \to 0 \ \ \textrm{as}\ \ n \to \infty
\end{align*}
for any $T>0$.
Part (3) then follows in the same way as in (3.18).

We now prove that $(u, v, \partial_t v) \in L_{loc}^{\infty} (\mathbb R ; D(\Delta) \oplus D(\Delta) \oplus H_0^1)$.
For any $t \in \mathbb R$ and $\psi \in D(\Delta)$, we have
\begin{align*}
   |(u(t) | \Delta \psi)| &= \lim_{n \to \infty} |(u_n(t) | \Delta \psi)|   \\
   &\leq \sup_{n} |(u_n(t) | \Delta \psi)|   \\
   &= \sup_{n} |(\Delta u_n(t) | \psi)|
   \leq M_2(t)^{1/2} \|\psi\|_2.
\end{align*}
This proves that $u(t) \in D(\Delta)$ and $\|\Delta u(t)\|_2^2 \leq M_2(t)$.
A similar argument shows that for any $t \in \mathbb R, (v(t), \partial_t v(t)) \in D(\Delta) \oplus H_0^1$ and $\|v(t)\|_{H^1}^2 + \|\partial_t v(t)\|_2^2 \leq M_2(t)$, as required.
This proves Part (4).

We next prove that $(u, v)$ is a unique global strong solution to (1.1)-(1.4) by showing that $(u, v)$ solves the associated integral equations
\begin{align*}
   &u(t) = U(t) \varphi + i \int_0^t U(t-s) v(s) u(s) ds,   \tag{5.2}   \\
   &v(t) = \dot{K}(t) \psi_0 + K(t) \psi_1 + \int_0^t K(t-s) |u(s)|^2 ds   \tag{5.3}
\end{align*}
as limits of (1.18) and (1.19), respectively.
We write
\begin{align*}
   &J_n^2 (J_n^2 v_n \cdot J_n^2 u_n) - v u   \\
   &= J_n^2 ( J_n^2 (v_n - v) \cdot J_n^2 u_n 
   + (J_n^2 v - v) \cdot J_n^2 u_n 
   + v J_n^2 (u_n - u)
   + v (J_n^2 u - u) )   \\
   &\hspace{11pt} + (J_n^2 (v u) - v u)   \tag{5.4}
\end{align*}
and estimate as
\begin{align*}
   &\|J_n^2 (J_n^2 v_n \cdot J_n^2 u_n) - v u\|_2   \\
   &\leq \|J_n^2 (v_n - v)\|_4 \|J_n^2 u_n\|_4 
   + \|J_n^2 v -v\|_4 \|J_n^2 u_n\|_4 
   + \|v\|_4 \|J_n^2 (u_n - u)\|_4   \\
   &\hspace{11pt} + \|v\|_4 \|J_n^2 u - u\|_4
   + \|J_n^2 (v u) - v u\|_2   \\
   &\leq C M_1^{1/2} ( \|v_n - v\|_{H^1} + \|J_n v - v\|_{H^1}
   + \|u_n - u\|_{H^1} + \|J_n u - u\|_{H^1} )   \\
   &\hspace{11pt} + 2 \|J_n (v u) - v u\|_2   \\
   &\leq C M_1^{1/2} ( \|v_n - v\|_{H^1} + \|u_n - u\|_{H^1} )
   + C n^{-1/2} M_1^{1/2} M_2^{1/2} + C n^{-1/2} M_2,   \tag{5.5}
\end{align*}
where we have used the estimates
\begin{align*}
   &\|J_n u - u\|_{H^1} \leq n^{-1/2} \|\nabla u\|_{H^1} \leq n^{-1/2} M_2^{1/2},   \\
   &\|J_n (v u) - v u\|_2 \leq n^{-1/2} \|\nabla (v u)\|_2 \leq n^{-1/2} M_2.
\end{align*}
By (5.4), we obtain
\begin{align*}
   &\sup_{|t| \leq T} \left\| \int_0^t U(t-s) J_n^2 (J_n^2 v_n(s) \cdot J_n^2 u_n(s)) ds
   - \int_0^t U(t-s) v(s) u(s) ds \right\|_2   \\
   &\leq \sup_{|t| \leq T} \left| \int_0^t \|J_n^2 (J_n^2 v_n(s) \cdot J_n^2 u_n(s)) 
   - v(s) u(s)\|_2 ds \right|   \\
   &\leq T \sup_{|t| \leq T} \|J_n^2 (J_n^2 v_n(t) \cdot J_n^2 u_n(t)) 
   - v(t) u(t)\|_2   \\
   &\leq C T M_1^{1/2} \sup_{|t| \leq T} (\|v_n(t) - v(t)\|_{H^1} + \|u_n(t) - u(t)\|_{H^1})
   + C n^{-1/2} M_2 \to 0   \tag{5.6}
\end{align*}
as $n \to \infty$ for any $T>0$.

Similarly, we write
\begin{align*}
   &J_n^2 |J_n^2 u_n|^2 - |u|^2   \\
   &= J_n^2 (\textrm{Re} (J_n^2 u_n - u) \overline{(J_n^2 u_n + u)})
   + (J_n^2 - I) |u|^2   \\
   &= J_n^2 (\textrm{Re} (J_n^2 (u_n - u) + (J_n^2 - I) u) \overline{(J_n^2 u_n + u)})
   + (J_n + I) (J_n - I) |u|^2   \\
   &= \textrm{Re} J_n^2 (\overline{(J_n^2 u_n - u)} (J_n^2 (u_n - u) + (J_n^2 - I) u))
   + (J_n + I) (J_n - I) |u|^2   \tag{5.7}
\end{align*}
and estimate as
\begin{align*}
   &\|J_n^2 |J_n^2 u_n|^2 - |u|^2\|_2   \\
   &\leq (\|J_n^2 u_n\|_4 + \|u\|_4) (\|J_n^2 (u_n - u)\|_4 + \|(J_n^2 - I) u\|_4)
   + 2 \|(J_n - I) |u|^2\|_2   \\
   &\leq C M_1^{1/2} (\|u_n - u\|_{H^1} + n^{-1/2} M_2^{1/2}) + C n^{-1/2} M_2. \tag{5.8}
\end{align*}
By (5.8), we obtain
\begin{align*}
   &\sup_{|t| \leq T} \left\| \int_0^t K(t-s) J_n^2 |J_n^2 u_n(s)|^2 ds
   - \int_0^t K(t-s) |u(s)|^2 ds \right\|_2   \\
   &\leq \sup_{|t| \leq T} \left| \int_0^t |t-s| \|J_n^2 |J_n^2 u_n(s)|^2 
   - |u(s)|^2\|_2 ds \right|   \\
   &\leq \frac{1}{2} T^2 \sup_{|t| \leq T} \|J_n^2 |J_n^2 u_n(s)|^2 
   - |u(s)|^2\|_2   \\
   &\leq C T^2 M_1^{1/2} \sup_{|t| \leq T} \|u_n(t) - u(t)\|_{H^1}
   + C T^2 n^{-1/2} M_2 \to 0   \tag{5.9}
\end{align*}
as $n \to \infty$ for any $T>0$.
Then (5.2) and (5.3) follow from (5.6) and (5.9), respectively.
It follows from Part (4) that the RHS of (5.2) belongs to $C(\mathbb R ; D(\Delta)) \cap C^1 (\mathbb R ; L^2)$ and that the RHS of (5.3) belongs to $C(\mathbb R ; D(\Delta)) \cap C^1 (\mathbb R ; H_0^1) \cap C^2 (\mathbb R ; L^2)$.
This means that the LHS of (5.2) and (5.3) belong to the corresponding classes, respectively, and therefore we see that $(u, v)$ is a global strong solutions to (1.1)-(1.4).
The uniqueness of global strong solutions follows by the standard Gronwall argument in $C( \mathbb R ; L^2)$.
This proves Part (1).   \\

\noindent
\textbf{Proof of Theorem 2.}   
The argument until (5.9) works also in the case $N=4$ since $H_0^1(\Omega) \hookrightarrow L^4(\Omega)$ holds for $N=4$.
Particularly, we have proved Parts (2)-(4) and the equalities (5.2) and (5.3).
Subsequent argument breaks down since we are not able to prove that $u v, |u|^2 \in L_{loc}^{\infty} (\mathbb R ; D(\Delta))$.
Instead, the weak continuity $(u, v, \partial_t v) \in C_w (\mathbb R ; D(\Delta) \oplus D(\Delta) \oplus H_0^1)$ follows from the strong continuity $(u, v, \partial_t v) \in C (\mathbb R ; H_0^1 \oplus H_0^1 \oplus L^2)$ and the local boundedness $(u, v, \partial_t v) \in L_{loc}^{\infty} (\mathbb R ; D(\Delta) \oplus D(\Delta) \oplus H_0^1)$.

We prove the uniqueness of solutions.
Let $(u, v, \partial_t v)$ and $(u', v', \partial_t v')$ be two triplets of solutions in $C_w (\mathbb R ; D(\Delta) \oplus D(\Delta) \oplus H_0^1)$ with the same data at $t=0$.
Let $T>0$ and $t \in [0,T]$.
We compute
\begin{align*}
   \frac{d}{dt} \|u - u'\|_2^2 &= 2 \textrm{Im} (i \partial_t (u - u') | u - u')   \\[-3pt]
   &= 2 \textrm{Im} (- \Delta (u - u') - (vu - v'u') | u - u')   \\
   &= - 2 \textrm{Im} ((v - v') u + v' (u - u') | u - u')   \\
   &= - 2 \textrm{Im} ((v - v') u | u - u')   \\
   &\leq 2 \|u\|_{\frac{2}{\varepsilon}} \|v - v'\|_{\frac{2}{1-\varepsilon /2}} \|u - u'\|_{\frac{2}{1-\varepsilon /2}}   \\
   &\leq C \varepsilon^{-1/2} \|u\|_{H^2} \|v - v'\|_{H^1}^{\varepsilon} \|v - v'\|_2^{1-\varepsilon} \|u - u'\|_{H^1}^{\varepsilon} \|u - u'\|_2^{1-\varepsilon}   \\
   &\leq C \varepsilon^{-1/2} M_2 M_1^{\varepsilon} \|v - v'\|_2^{1-\varepsilon} \|u - u'\|_2^{1-\varepsilon}   \\
   &\leq C \varepsilon^{-1/2} M_2 (1 + M_1) (\|u - u'\|_2^{2-2\varepsilon} + \|v - v'\|_2^{2-2\varepsilon})   \\
   &\leq C \varepsilon^{-1/2} M_2 (1 + M_1) (\|u - u'\|_2^2 + \|v - v'\|_2^2)^{1-\varepsilon}   \tag{5.10}
\end{align*}
for any $\varepsilon \in (0, \frac{1}{2}]$, where $C$ is independent of $\varepsilon$ and we have used the H\"{o}lder inequality with $1 = \frac{\varepsilon}{2} + \frac{1-\varepsilon /2}{2} + \frac{1-\varepsilon /2}{2}$, the Sobolev type inequality \cite{Og,OO,Oz}
\begin{align*}
   \|u\|_p \leq C \sqrt{p} \|u\|_{H^2}   \tag{5.11}
\end{align*}
for any $p \in [2, + \infty)$ with $p = \frac{2}{\varepsilon}$, Lemma 2 with $\delta_4(\frac{2}{1-\varepsilon /2}) = \varepsilon$, and $M_2 := M_2(T)$ for $N=4$.
As regards $(v, \partial_t v)$ and $(v', \partial_t v')$, we compute
\begin{align*}
   &\frac{d}{dt} (\|v - v'\|_2^2 + \|\omega^{-1} \partial_t (v - v')\|_2^2)   \\
   &= 2 (\omega^{-1} (|u|^2 - |u'|^2) | \omega^{-1} \partial_t (v - v'))   \\
   &= 2 \textrm{Re} (\omega^{-1} ((\overline{u} + \overline{u'}) (u - u')) | \omega^{-1} \partial_t (v - v'))   \\
   &\leq 2 \|\omega^{-1} ((\overline{u} + \overline{u'}) (u - u'))\|_2 \|\omega^{-1} \partial_t (v - v')\|_2   \\
   &\leq C \|(\overline{u} + \overline{u'}) (u - u')\|_{4/3} \|\omega^{-1} \partial_t (v - v')\|_2   \\
   &\leq C (\|u\|_4 + \|u'\|_4) \|u - u'\|_2 \|\omega^{-1} \partial_t (v - v')\|_2   \\
   &\leq C M_1^{1/2} (\|u - u'\|_2^2 + \|\omega^{-1} \partial_t (v - v')\|_2^2),   \tag{5.12}
\end{align*}
where we have used the embedding $L^{4/3} \hookrightarrow H^{-1}$ as the dual of the embedding $H_0^1 \hookrightarrow L^4$ for $N=4$.
We define
\begin{align*}
   F(t) := \|u(t) - u'(t)\|_2^2 + \|v(t) - v'(t)\|_2^2 + \|\omega^{-1} \partial_t (v - v')(t)\|_2^2.
\end{align*}
By (5.10) and (5.12), $F$ satisfies the differential inequality
\begin{align*}
   F'(t) \leq C \varepsilon^{1/2} F(t)^{1-\varepsilon} + C F(t),   \tag{5.13}
\end{align*}
where we have omitted the dependence of constants on $M_1$ and $M_2$.
With any $\delta \in (0, 1]$, we obtain from (5.13)
\begin{align*}
   ((F + \delta)^{\varepsilon})'(t) &= \varepsilon (F(t) + \delta)^{\varepsilon -1} F'(t)   \\
   &\leq C \varepsilon^{1/2} + C \varepsilon (F(t) + \delta)^{\varepsilon}.   \tag{5.14}
\end{align*}
This yields
\begin{align*}
   (\exp (-C \varepsilon t) (F + \delta)^{\varepsilon})'(t) 
   \leq \exp (-C \varepsilon t) C \varepsilon^{1/2} 
   \leq C \varepsilon^{1/2}.   \tag{5.15}
\end{align*}
Integrating both sides of (5.15), we have
\begin{align*}
   \exp (-C \varepsilon t) (F(t) + \delta)^{\varepsilon} 
   \leq \delta^{\varepsilon} + C \varepsilon^{1/2} T,   \tag{5.16}
\end{align*}
since $F(0) = 0$.
By (5.16) and letting $\delta \downarrow 0$, we obtain
\begin{align*}
   \sup_{t \in [0,T]} F(t) \leq \exp (CT) (C \varepsilon^{1/2} T)^{1/2}.   \tag{5.17}
\end{align*}
For any $T>0$, there exists $\varepsilon_0 \in (0, \frac{1}{2}]$ such that $C \varepsilon_0^{1/2} T \leq \frac{1}{2}$.
Then letting $\varepsilon \downarrow 0$ on the RHS of (5.17), we obtain the uniqueness.
This completes the proof of Theorem 2.

%%%%% 6. Remarks on global finite energy solutions %%%%%

\section{Remarks on global finite energy solutions}

\ \ \ In this section we give some remarks on the finite energy solutions, namely solutions $(u, v) \in C(\mathbb R ; H_0^1 \oplus H_0^1)$ of the integral equations (5.2)-(5.3) with conserved charge and energy.
As in the case of strong solutions, the situation depends essentially on the space dimension $N$ via Sobolev embedding.
We start with the simplest case $N=1$ where the embedding $H^1(\Omega) \hookrightarrow L^{\infty}(\Omega)$ holds.   \\

\noindent
\textbf{Theorem 3.}\ \ \textit{Let $N=1$ and let $(\varphi, \psi_0, \psi_1) \in H_0^1 \oplus H_0^1 \oplus L^2$.
Then}:   \\
(1) \textit{There exists a unique pair of solutions $(u, v)$ satisfying} (5.2)-(5.3),
\begin{align*}
   &u \in C(\mathbb R ; H_0^1) \cap C^1(\mathbb R ; H^{-1}),   \\
   &v \in C(\mathbb R ; H_0^1) \cap C^1(\mathbb R ; L^2) \cap C^2(\mathbb R ; H^{-1}),
\end{align*}
\textit{where $H^{-1} = H^{-1}(\Omega) = (H_0^1(\Omega))'$.}   \\
(2) \textit{The total charge $Q(t)$ and energy $E(t)$ are conserved in time.}   \\
(3) \textit{The estimate} (1.8) \textit{holds.}   \\

\noindent
\textbf{Remark 2.}\ \ Tsutsumi and Fukuda \cite{FT1} proved Theorem 3 in the case $\Omega = \mathbb R$.
The method of proof in \cite{FT1} depends on a compactness argument of Galerkin type, while the proof below is independent of compactness arguments.   \\

\noindent
\textbf{Proof of Theorem 3.}\ \ Regarding the sequence $((u_n, v_n, \partial_t v_n) ; n \geq 1)$ of solutions of (1.18)-(1.19) with $(\varphi, \psi_0, \psi_1) \in H_0^1 \oplus H_0^1 \oplus L^2$, the argument on the boundedness in $C(\mathbb R ; H_0^1 \oplus H_0^1 \oplus L^2)$ and the convergence in $C(\mathbb R ; L^2 \oplus L^2 \oplus H^{-1})$ in Section 3 and 4 holds in the case $N=1$.
We denote by $(u, v, w) \in C(\mathbb R ; L^2 \oplus L^2 \oplus H^{-1})$ its limit.
By the same argument as in Section 5, we have $w= \partial_t v$.
By the Sobolev embedding $H^1 \hookrightarrow L^{\infty}$, we see that $J_n^2 (J_n^2 v_n \cdot J_n^2 u_n)$ and $J_n^2 |J_n^2 u_n|^2$ tend to $v u$ and $|u|^2$ in $C(\mathbb R ; L^2)$, respectively.
This proves the existence of solutions to (5.2)-(5.3).
The uniqueness follows from the standard Gronwall argument.
This proves Part (1).
A similar argument shows that $(|J_n^2 u_n|^2 | J_n^2 v_n)$ tends to $(|u|^2 | v)$ in $C(\mathbb R ; \mathbb R)$.
This proves Part (2).
We have already proved (1.8) in the setting of Theorem 3.
This proves Part (3).   \\

We next study the case $N=2$, where the embedding $H^1(\Omega) \hookrightarrow L^{\infty}(\Omega)$ is no longer available, while there exists $C>0$ such that for any $p \in [2, +\infty)$ and any $u \in H_0^1(\Omega)$ the estimate
\begin{align*}
   \|u\|_p \leq C p^{1/2} \|u\|_{H^1}   \tag{6.1}
\end{align*}
holds \cite{Og,OO,Oz}.   \\

\noindent
\textbf{Theorem 4.}\ \ \textit{Let $N=2$ and let $(\varphi, \psi_0, \psi_1) \in H_0^1 \oplus H_0^1 \oplus L^2$.
Then, all the statements in Theorem 3 hold.}   \\

\noindent
\textbf{Remark 3.}\ \ The method of proof depends on the estimate (6.1).
Similar methods may be found in \cite{FMO,Hm,OO,OT1,OT2}.   \\

\noindent
\textbf{Proof of Theorem 4.}\ \ For $n \in \mathbb Z_{>0}$ we consider the Cauchy problem
\begin{align*}
   &i \partial_t u + \Delta u = - v u,   \tag{6.2}   \\
   &\partial_t^2 v - \Delta v + v = |u|^2,   \tag{6.3}   \\
   &(u(0), v(0), \partial_t v(0)) = (J_n \varphi, J_n \psi_0, J_n \psi_1).   \tag{6.4}
\end{align*}
Since $(J_n \varphi, J_n \psi_0, J_n \psi_1) \in D(\Delta) \oplus D(\Delta) \oplus H_0^1$, Theorem 1 shows the existence and uniqueness of global strong solutions to (6.2)-(6.4) with conserved charge and energy.
We denote by $(u_n, v_n, \partial_t v_n)$ the corresponding solution with index $n$.
By a similar and simpler argument in Section 2, we see that the sequence $((u_n, v_n, \partial_t v_n) ; n \geq 1)$ is bounded in $L^{\infty} (\mathbb R ; H_0^1 \oplus H_0^1 \oplus L^2)$ with uniform bound
\begin{align*}
   \sup_{n \geq 1} \sup_{t \in \mathbb R} 
   (\|u_n(t)\|_{H^1}^2 + \|v_n(t)\|_{H^1}^2 + \|\partial_t v_n(t)\|_2^2)
   \leq M_1 < \infty.
\end{align*}
We now prove that $((u_n, v_n, \partial_t v_n) ; n \geq 1)$ is a convergent sequence in $C([-T, T] ; L^2 \oplus L^2 \oplus H^{-1})$ for any $T>0$.
For that purpose, with $m>n$ we compute
\begin{align*}
  \frac{d}{dt} \|u_m - u_n\|_2^2
  &= 2 \textrm{Im} (i \partial_t (u_m - u_n) | u_m - u_n)   \\
  &= 2 \textrm{Im} (- \Delta (u_m - u_n) - (v_m u_m - v_n u_n) | u_m - u_n)   \\
  &= - 2 \textrm{Im} ((v_m - v_n) u_m + v_n (u_m - u_n) | u_m - u_n)   \\
  &= - 2 \textrm{Im} ((v_m - v_n) u_m | u_m - u_n)   \\
  &\leq 2 \|u_m\|_{\frac{1}{\varepsilon}} \|v_m - v_n\|_{\frac{2}{1-\varepsilon}} 
  \|u_m - u_n\|_{\frac{2}{1-\varepsilon}}   \\
  &\leq C \varepsilon^{-\frac{1}{2}} \|u_m\|_{H^1} \|v_m - v_n\|_{H^1}^{\varepsilon} 
  \|v_m - v_n\|_2^{1-\varepsilon} \|u_m - u_n\|_{H^1}^{\varepsilon} 
  \|u_m - u_n\|_2^{1-\varepsilon}   \\
  &\leq C \varepsilon^{-\frac{1}{2}} M_1^{\frac{1}{2} + \varepsilon}
  \|v_m - v_n\|_2^{1-\varepsilon} \|u_m - u_n\|_2^{1-\varepsilon}   \\
  &\leq C \varepsilon^{-\frac{1}{2}} (M_1^{\frac{1}{2}} + M_1)
  (\|v_m - v_n\|_2^{2-2\varepsilon} + \|u_m - u_n\|_2^{2-2\varepsilon})   \\
  &\leq C \varepsilon^{-\frac{1}{2}} (M_1^{\frac{1}{2}} + M_1)
  (\|v_m - v_n\|_2^2 + \|u_m - u_n\|_2^2)^{1-\varepsilon},   \tag{6.5}
\end{align*}
where $C>0$ is independent of $m,n$, and $\varepsilon \in (0, \frac{1}{2}]$ and we have used the H\"{o}lder inequality with $1 = \varepsilon + \frac{1-\varepsilon}{2} + \frac{1-\varepsilon}{2}$, (6.1) with $p=\frac{1}{\varepsilon}$, and Lemma 2 with $\delta_2 (\frac{1-\varepsilon}{2}) = \varepsilon$.
As regards $(v_n, \partial_t v_n)$, we compute
\begin{align*}
   &\frac{d}{dt} (\|v_m - v_n\|_2^2 + \|\omega^{-1} \partial_t (v_m - v_n)\|_2^2)   \\
   &= 2 (\omega^{-1} (|u_m|^2 - |u_n|^2) | \omega^{-1} \partial_t (v_m - v_n))   \\
   &= 2 \textrm{Re} (\omega^{-1} ((\overline{u_m} + \overline{u_n}) (u_m - u_n)) 
   | \omega^{-1} \partial_t (v_m - v_n))   \\
   &\leq 2 \|\omega^{-1} ((\overline{u_m} + \overline{u_n}) (u_m - u_n))\|_2
   \|\omega^{-1} \partial_t (v_m - v_n)\|_2   \\
   &\leq C \|(\overline{u_m} + \overline{u_n}) (u_m - u_n)\|_{4/3}
   \|\omega^{-1} \partial_t (v_m - v_n)\|_2   \\
   &\leq C (\|u_m\|_4 + \|u_n\|_4) \|u_m - u_n\|_2
   \|\omega^{-1} \partial_t (v_m - v_n)\|_2   \\
   &\leq C M_1^{1/2} (\|u_m - u_n\|_2^2 +
   \|\omega^{-1} \partial_t (v_m - v_n)\|_2^2),   \tag{6.6}   
\end{align*}
where we have used the embedding $L^{\frac{4}{3}} \hookrightarrow H^{-1}$ as the dual of the embedding $H_0^1 \hookrightarrow L^4$.
We define
\begin{align*}
   F_{m,n}(t) := \|u_m(t) - u_n(t)\|_2^2 + \|v_m(t) - v_n(t)\|_2^2
   + \|\omega^{-1}  \partial_t (v_m - v_n)(t)\|_2^2.
\end{align*}
By (6.5) and (6.6), $F_{m,n}$ satisfies the differential inequality
\begin{align*}
   F_{m,n}'(t) \leq C \varepsilon^{-1/2} F_{m,n}(t)^{1-\varepsilon} + C F_{m,n}(t),   \tag{6.7}
\end{align*}
where we have omitted the dependence of constants on $M_1$.
With any $\delta \in (0,1]$, we obtain from (6.7)
\begin{align*}
   ((F_{m,n} + \delta)^{\varepsilon})'(t) 
   &= \varepsilon (F_{m,n}(t) + \delta)^{\varepsilon -1} F_{m,n}'(t)   \\
   &\leq C \varepsilon^{1/2} + C \varepsilon (F_{m,n}(t) + \delta)^{\varepsilon}.   \tag{6.8}
\end{align*}
This yields
\begin{align*}
   (\exp (-C \varepsilon t) (F_{m,n} + \delta)^{\varepsilon})'(t)
   \leq \exp (-C \varepsilon t) C \varepsilon^{1/2} \leq C \varepsilon^{1/2}.   \tag{6.9}
\end{align*}
Let $T>0$ and $t \in [0,T]$.
Then integrating both sides of (6.9), we have
\begin{align*}
   \exp (-C \varepsilon t) (F_{m,n}(t) + \delta)^{\varepsilon}
   \leq (F_{m,n}(0) + \delta)^{\varepsilon} + C \varepsilon^{1/2} T.   \tag{6.10}
\end{align*}
By (6.10) and letting $\delta \downarrow 0$, we obtain
\begin{align*}
   \sup_{t \in [0,T]} F_{m,n}(t) 
   \leq \exp (CT) (F_{m,n}(0)^{\varepsilon} + C \varepsilon^{1/2} T)^{1/\varepsilon}.   \tag{6.11}
\end{align*}
In the same way as in the argument in \cite{OT1}, (6.11) implies
\begin{align*}
   \sup_{t \in [0,T]} F_{m,n}(t) 
   \leq \exp (CT) \left( \frac{1}{2} + C \varepsilon^{1/2} T \right)^{1/\varepsilon},   \tag{6.12}
\end{align*}
for $m, n$ large enough.
The RHS of (6.12) tends to zero as $\varepsilon \downarrow 0$.

This proves that the sequence $((u_n, v_n, \partial_t v_n) ; n \geq 1)$ has a limit $(u, v, w) \in C(\mathbb R ; L^2 \oplus L^2 \oplus H^{-1})$.
In the same way as above, we have $\partial_t v = w$.
Moreover, for any $T>0$
\begin{align*}
   \sup_{t \in [0,T]} \|v_n u_n - v u\|_{H^{-1}} 
   &\leq C \sup_{t \in [0,T]} \| v_n u_n - v u\|_{4/3}   \\
   &\leq C \sup_{t \in [0,T]} (\|v_n - v\|_2 \|u_n\|_4 + \|v\|_4 \|u_n - u\|_2)   \\
   &\leq C M_1^{1/2} \sup_{t \in [0,T]} (\|v_n - v\|_2 + \|u_n - u\|_2) \to 0, \tag{6.13} \\
   \sup_{t \in [0,T]} \||u_n|^2 - |u|^2\|_{H^1}
   &\leq C M_1^{1/2} \sup_{t \in [0,T]} \|u_m - u_n\|_2 \to 0   \tag{6.14}
\end{align*}
as $n \to \infty$.
By (6.13) and (6.14), we find that $(u,v)$ satisfies (5.2) and (5.3) in $C(\mathbb R ; H^{-1})$.
The rest of the statements in Theorem 4 follow in the same way as before.

%%%%%%%%%%%%%%%%%%%% References %%%%%%%%%%%%%%%%%%%%

%%%%%%%%%%%%%%%%%%%%%%%%%% 参考 %%%%%%%%%%%%%%%%%%%%%%%%%%%%%%%
%    LATEX2イプシロン 美文書作成入門, 奥村晴彦, 黒木裕介, 技術評論社     
%    LATEX2イプシロン 辞典, 吉永徹美, 翔泳社                                       
%    TEX Wiki https://texwiki.texjp.org
%%%%%%%%%%%%%%%%%%%%%%%%%%%%%%%%%%%%%%%%%%%%%%%%%%%%%%%%%%%%%%

\end{document}